\documentclass[12pt,leqno]{amsart}
\usepackage{amsmath,amssymb,amsthm,amsfonts}

\newtheorem{thm}{Theorem}[section]
\newtheorem{prop}[thm]{Proposition}
\newtheorem{coro}[thm]{Corollary}
\newtheorem{lem}[thm]{Lemma}

\theoremstyle{definition}
\newtheorem{defi}[thm]{Definition}
\newtheorem{ex}[thm]{Example}

\theoremstyle{remark}
\newtheorem{rem}[thm]{Remark}

\numberwithin{equation}{section}

\newcommand{\nc}{\newcommand}

\nc{\op}{\oplus}
\nc{\ot}{\otimes}
\nc{\bop}{\bigoplus}
\nc{\ze}{\{ 0 \}}
\nc{\pv}{P^{\vee}}

\nc{\J}{\mathbf{J}}
\nc{\Q}{\mathbf{Q}}
\nc{\Z}{\mathbf{Z}}
\nc{\m}{\mathbf{m}}

\nc{\g}{{\mathfrak{g}}}
\nc{\h}{\mathfrak{h}}

\nc{\fit}{\tilde f_i}
\nc{\eit}{\tilde e_i}

\nc{\vep}{\varepsilon}
\nc{\vp}{\varphi}

\nc{\cl}{\colon}
\nc{\re}{\mathrm{re}}
\nc{\im}{\mathrm{im}}

\nc{\Uq}{U_q (\g)}
\nc{\Um}{{U_q^-(\g)}}
\nc{\CO}{\mathcal{O}}

\nc{\vla}{V(\la)}
\nc{\lla}{L(\la)}
\nc{\bla}{B(\la)}
\nc{\gla}{G(\la)}

\nc{\vmu}{V(\mu)}
\nc{\lmu}{L(\mu)}
\nc{\bmu}{B(\mu)}

\nc{\ltau}{L(\tau)}

\nc{\vlm}{V(\la+\mu)}
\nc{\llm}{L(\la+\mu)}
\nc{\blm}{B(\la+\mu)}

\nc{\lom}{\lla \ot \lmu}

\nc{\al}{\alpha}
\nc{\be}{\beta}
\nc{\ga}{\gamma}
\nc{\del}{\delta}
\nc{\Del}{\Delta}
\nc{\La}{\Lambda}
\nc{\la}{\lambda}

\nc{\ophi}{\Phi_{\la,\mu}}
\nc{\opsi}{\Psi_{\la,\mu}}
\nc{\os}{S_{\la,\mu}}

\nc{\add}[3]{\sum_{#1=#2}^{#3}}
\nc{\wf}[1]{\tilde f_{i_{#1}}}
\nc{\bbin}[2]
{{\left\{ \begin{matrix} #1 \\ #2 \end{matrix} \right\}}}
\nc{\isoto}[1][]{\xrightarrow[#1]{\sim}}
\nc{\str}[1][i]{$#1$-string decomposition}

\nc{\A} {{\mathbf A}}
\nc{\Az}{{\mathbf A_0}}
\nc{\Ai}{{\mathbf A_\infty}}
\nc{\Li}{L_\infty}

\nc{\linf}{L(\infty)}
\nc{\binf}{B(\infty)}
\nc{\ginf}{G(\infty)}

\nc{\ovl}{\overline}
\nc{\mt}{\mapsto}

\nc{\lan}{\langle}
\nc{\ran}{\rangle}

\nc{\seteq}{\operatorname{{:=}}}
\nc{\set}[2]{\left\{#1\,;\,#2\right\}}

\renewcommand{\labelenumii}{(\roman{enumii})}

\newenvironment{tenumerate}
{\begin{enumerate}

}{\end{enumerate}}

\renewcommand{\ker}{\operatorname{Ker}}

\nc{\wt}{\operatorname{wt}}
\nc{\id}{\operatorname{id}}

\begin{document}

\title[Crystal Bases for Quantum Generalized Kac-Moody Algebras]
      {Crystal Bases for Quantum Generalized \\ Kac-Moody Algebras}

\author[K. Jeong, S.-J. Kang, M. Kashiwara]
{Kyeonghoon Jeong$^{\diamond}$,
Seok-Jin Kang$^{*}$,
Masaki Kashiwara$^{\dagger}$}

\address {\parbox{400pt}{$^{\diamond}$Department of Mathematics, NS30,
Seoul National University,
Seoul 151-747, Korea\\}}

\email{khjeong@math.snu.ac.kr}

\address{\parbox{400pt}{$^{*}$School of Mathematics,
Korea Institute for Advanced Study,
207-43 Cheongryangri-Dong,
Dongdaemun-Gu,
Seoul 130-722, Korea\\}}

\email{sjkang@kias.re.kr}

\address{\parbox{400pt}{$^{\dagger}$Research Institute for Mathematical Sciences,
Kyoto University, Kyoto 606, Japan\\}}

\email{masaki@kurims.kyoto-u.ac.jp}

\thanks{$^{\diamond}$This research was supported by KOSEF Grant
\# 98-0701-01-5-L}

\thanks{$^{*}$This research was supported by KOSEF Grant
\# 98-0701-01-5-L and the Young Scientist Award,
Korean Academy of Science and Technology}

\thanks{$^{\dagger}$This research is partially supported by
Grant-in-Aid for Scientific Research (B1)13440006,
Japan Society for the Promotion of Science.}

\begin{abstract}
In this paper, we develop the crystal basis theory for
quantum generalized Kac-Moody algebras.
For a quantum generalized Kac-Moody algebra $\Uq$,
we first introduce the category $\CO_{int}$ of $\Uq$-modules
and prove its semisimplicity.
Next, we define the notion of crystal bases for $\Uq$-modules
in the category $\CO_{int}$ and for the subalgebra $\Um$.
We then prove the tensor product rule and
the existence theorem for crystal bases.
Finally, we construct the global bases for $\Uq$-modules
in the category $\CO_{int}$ and for the subalgebra $\Um$.
\end{abstract}

\maketitle

\section*{{\bf Introduction}}

The {\it quantum groups}, introduced independently by Drinfel'd and Jimbo,
are certain families of Hopf algebras that are
deformations of universal enveloping algebras of
Kac-Moody algebras \cite{Dr85, J85}.
More precisely, let $\g$ be a symmetrizable Kac-Moody algebra
and let $U(\g)$ be its universal enveloping algebra.
Then, for each generic parameter $q$, we can associate a
Hopf algebra $\Uq$, called the quantum group,
whose structure {\it tends to} that of $U(\g)$ as $q$ approaches 1.
In \cite{Lus88}, Lusztig showed that
the integrable highest weight modules over $U(\g)$ can be
deformed to those over $\Uq$ in such a way that the dimensions
of weight spaces are invariant under the deformation.

In \cite{Kas90, Kas91}, Kashiwara
developed the {\it crystal basis theory} for the quantum groups
associated with symmetrizable Kac-Moody algebras.
(In \cite{Lus90}, a more geometric approach was developed by
Lusztig, which is called the {\it canonical basis theory}.)
The {\it crystal bases} can be understood as bases at $q=0$ and they are
given a structure of colored oriented graphs, called the
{\it crystal graphs}.
The crystal graphs have many nice combinatorial features
reflecting the internal structure of integrable modules
over quantum groups.
In particular, they have extremely simple behavior with respect to
taking the tensor product.
Thus the crystal basis theory provides us with a very powerful
combinatorial tool to investigate the structure of integrable
modules over quantum groups.

On the other hand, in his study of Monstrous moonshine,
Borcherds introduced a new class of infinite
dimensional Lie algebras called the {\it generalized Kac-Moody algebras}
\cite{B88, B92}.
The structure and the representation theory of generalized Kac-Moody
algebras are very similar to those of Kac-Moody algebras, and many
basic facts about Kac-Moody algebras can be extended to generalized
Kac-Moody algebras.
But there are some differences, too.
For example, generalized Kac-Moody algebras may have
{\it imaginary simple roots}
with norms $\le 0$ whose multiplicity
can be greater than one.
Also, they may have infinitely many simple roots.

In \cite{Ka95}, Kang constructed the quantum groups $\Uq$
associated with generalized Kac-Moody algebras $\g$.
For simplicity, we will call $\Uq$ the {\it quantum
generalized Kac-Moody algebras}.
Moreover, he showed that, for a generic $q$, Verma modules and
irreducible highest weight modules over $U(\g)$ with dominant integral
highest weights can be deformed to those over $\Uq$ in such a way that
the dimensions of weight spaces are invariant under the deformation.
These results were extended to generalized Kac-Moody superalgebras
by Benkart, Kang and Melville \cite{BKM98}.

In this paper, we develop the crystal basis theory for quantum
generalized Kac-Moody algebras following the framework given in
\cite{Kas91}.
After reviewing some of the basic facts about quantum generalized
Kac-Moody algebras and their modules,
we define the category $\CO_{int}$ consisting of
integrable $\Uq$-modules satisfying certain conditions
on their weights and show that this category is semisimple.
We then define the notion of crystal bases for $\Uq$-modules
in the category $\CO_{int}$ and prove the standard properties
of crystal bases including the {\it tensor product rule}.

Next, we prove the existence theorem for crystal bases using
Kashiwara's {\it grand-loop argument}, which consists of 15
interlocking inductive statements.
As is the case with Kac-Moody algebras,
the fundamental properties of crystal bases for $\Um$
play a crucial role in this argument.

Finally, we {\it globalize} the theory of crystal bases.
That is, we show that any irreducible highest weight
$\Uq$-module in the category $\CO_{int}$ has a unique {\it global basis}.
The existence of global bases is proved
using the basic properties of {\it balanced triple}
and the triviality of vector bundles over $\mathbf{P}^{1}$.

\vskip 5mm

\noindent
{\bf Acknowledgments.}  We would like to express our sincere gratitude
to Olivier Schiffmann for many valuable discussions.

\medskip

\section{{\bf Quantum Generalized Kac-Moody Algebras}}

Let $I$ be a finite or countably infinite index set.
A real square matrix $A=(a_{ij})_{i,j \in I}$ is called a
{\it Borcherds-Cartan matrix} if it satisfies:
\begin{tenumerate}
\item $a_{ii}=2$ or $a_{ii}\le 0$ for all $i \in I$,
\item $a_{ij} \le 0$ if $i \ne j$,
\item $a_{ij}\in \Z$ if $a_{ii}=2$,
\item $a_{ij}=0$ if and only if $a_{ji} = 0$.
\end{tenumerate}
We say that an index $i$ is {\it real} if $a_{ii}=2$ and
{\it imaginary} if $a_{ii}\le 0$.
We denote by $I^\re \seteq \set{i \in I}{a_{ii}=2}$,
$I^\im \seteq \set{i \in I}{a_{ii} \le 0}$.

In this paper, we assume that all the entries of $A$ are
integers and the diagonal entries are even.
Furthermore, we assume that $A$ is {\it symmetrizable}; that is,
there is a diagonal matrix
$D= \mbox{diag} (s_i \in \Z_{>0} \mid i \in I)$
such that $DA$ is symmetric.

Let $\pv$ be a free abelian group
and let $\Pi^\vee = \{ h_i \}_{i \in I}$ be a family of elements
in $\pv$ called the set of {\it simple coroots},
and set $\h = \Q \ot_{\Z} \pv$.
The free abelian group $\pv$ is called the
{\it co-weight lattice} and the $\Q$-vector space $\h$
is called the {\it Cartan subalgebra}.
The {\it weight lattice} is defined to be
\begin{align*}
P \seteq \set{\la \in \h^*}{\la(\pv) \subset \Z}.
\end{align*}
We assume that $P$ contains
a family of elements $\Pi = \{ \al_i \}_{i \in I }$ called
the set of {\it simple roots} satisfying
\begin{align*}
\al_i(h_j) = a_{ji} \quad \mbox{for } i, j \in I.
\end{align*}

\begin{defi}
We call the quintuple $(A, \pv, P, \Pi^\vee, \Pi)$
a {\it Borcherds-Cartan datum} associated with $A$.
\end{defi}

\smallskip

In this paper, we assume that
there is a non-degenerate symmetric bilinear form $(\ | \ )$
on $\h^*$ satisfying
\begin{align*}
(\al_i | \la) = s_i \la (h_i)\quad\mbox{for every } \la \in \h^*,
\end{align*}
and therefore we have
\begin{align*}
(\al_i | \al_j) = s_i a_{ij} \quad \mbox{for every } i, j \in I.
\end{align*}
We also assume that there exists $\La_i \in P$ such that
$\La_i(h_j)=\del_{ij}$ for $i, j \in I$.
Such a $\La_i$ is called a {\em fundamental weight}.
We assume further that $\{ \al_i \}_{i \in I}$ is linearly independent.

\smallskip

For a symmetrizable Borcherds-Cartan matrix, there
exists such a Borcherds-Cartan datum.
Take
\begin{align*}
P^\vee= \left( \bop_{i \in I} \Z h_i \right)
\op \left( \bop_{i \in I} \Z d_i \right)
\end{align*}
and define
$\al_i$ and $\La_i$ by
\begin{align*}
&\al_i(h_j)=a_{ji} \mbox{ and }\al_i(d_j)=\del_{ij},\\
&\La_i(h_j)=\del_{ij} \mbox{ and } \La_i(d_j)=0.
\end{align*}
Let $(\ | \ )$ be the bilinear form on
$\bigl(\op_i(\Q \al_i\op \Q\La_i)\bigr)\times\h^*$ defined by
\begin{align*}
(\al_i | \la) = s_i \la(h_i), \quad (\La_i | \la) = s_i \la(d_i).
\end{align*}
Since it is symmetric on
$\bigl(\op_i(\Q \al_i\op \Q\La_i)\bigr)
\times\bigl(\op_i(\Q \al_i\op \Q \La_i)\bigr)$,
we can extend this to a symmetric form on $\h^*$.
Then such a symmetric form is non-degenerate.

\smallskip

We denote by $P^+$
the set $\set{\la \in P}{\la(h_i) \ge 0 \mbox{ for every }i \in I}$
of {\it dominant integral weights}.
The free abelian group $Q \seteq \bop_{i \in I} \Z \al_i$ is called the
{\it root lattice}.
We set $Q_+ = \sum_{i \in I} \Z_{\ge 0} \al_i$ and
$Q_{-} = -Q_+$.
For $\al \in Q_+$, we can write $\al=\sum_{k=1}^r\al_{i_k}$
for $i_1,\dotsc,i_r\in I$.
We set $|\al|=r$ and call it the {\it height} of $\al$.

For an indeterminate $q$, set $q_i = q^{s_i}$ and define
\begin{align}
[n]_i = \dfrac{q_i^{n} - q_i^{-n}} {q_i - q_i^{-1}},
\quad \quad [n]_i! = \prod_{k=1}^n [k]_i, \quad \quad
\left[\begin{matrix} m \\ n \end{matrix}\right]_i
= \dfrac{[m]_i!}{[m-n]_i! [n]_i!}.
\end{align}
If $a_{ii}<0$, set $c_i = - \dfrac{1}{2} a_{ii} \in \Z_{>0}$,
and define
\begin{align}
\{ n \}_i = \dfrac{q_i^{c_i n} - q_i^{-c_i n}}
{q_i^{c_i} - q_i^{-c_i}},
\quad \quad \{ n \}_i! = \prod_{k=1}^n \{ k \}_i, \quad \quad
{\bbin m n}_i
= \dfrac{\{ m \}_i!}{\{ m-n \}_i! \{ n \}_i!}.
\end{align}
For convenience, if $a_{ii}=0$, we will also use the notation
$c_i = 0$ and
\begin{align}
\{ n \}_i = n,
\quad \quad \{ n \}_i! = n !, \quad \quad
\bbin m n _i
= \binom m n .
\end{align}
\begin{defi}
The {\it quantum generalized Kac-Moody algebra} $\Uq$
associated with a Borcherds-Cartan datum
$(A, \pv, P, \Pi^\vee, \Pi)$ is the associative algebra
over $\Q(q)$ with the unit $1$
generated by the symbols $e_i$, $f_i$ $(i \in I)$ and $q^h$
$(h\in \pv)$ subject to the following defining relations\,:

{\allowdisplaybreaks
\begin{align}
\begin{aligned}
\ & q^0 = 1, \ \
q^h q^{h'} = q^{h+ h'} \quad (h, h'\in \pv),\\
\ & q^h e_i q^{-h} = q^{\al_i(h)} e_i, \quad
q^h f_i q^{-h} = q^{-\al_i(h)} f_i, \\
\ & e_i f_j - f_j e_i
= \del_{ij} \frac{K_i - K_i^{-1}}{q_i - q_i^{-1}},
\quad \mbox{where} \quad K_i = q^{s_i h_i}, \\
\ & \sum_{r=0}^{1-a_{ij}} (-1)^r
{\begin{bmatrix} 1-a_{ij} \\ r \end{bmatrix}}_i
e_i^{1-a_{ij}-r} e_j e_i^{r} = 0 \quad \mbox{if}
\ \ a_{ii}=2, \ i \ne j, \\
\ & \sum_{r=0}^{1-a_{ij}} (-1)^r
{\begin{bmatrix} 1-a_{ij} \\ r \end{bmatrix}}_i
f_i^{1-a_{ij}-r} f_j f_i^{r} = 0 \quad \mbox{if}
\ \ a_{ii}=2, \ i \ne j, \\
\ & e_ie_j-e_je_i=0, \quad
f_if_j - f_j f_i = 0
\quad \mbox{if} \ \ a_{ij}=0.
\end{aligned}
\end{align}
}
\end{defi}
From now on, we will use the notation
\begin{align}
\begin{aligned}
& e_i^{(n)} = \dfrac{e_i^n}{[n]_i!},
\quad f_i^{(n)} = \dfrac{f_i^n}{[n]_i!}
\quad \mbox{if} \ \ a_{ii} = 2, \\
& e_i^{(n)} = e_i^n,
\quad f_i^{(n)} = f_i^n
\quad \mbox{if} \ \ a_{ii} \le 0.
\end{aligned}
\end{align}
We understand that
$ e_i^{(n)}=f_i^{(n)}=1$ for $n=0$ and $ e_i^{(n)}=f_i^{(n)}=0$ for $n<0$.
For $\al = \sum k_i \al_i\in Q$, we will
also use the notation
\begin{align*}
K_{\al} = \prod K_i^{k_i}.
\end{align*}

The quantum generalized Kac-Moody algebra
$\Uq$ has a {\it Hopf algebra} structure
with the comultiplication $\Del$, the counit $\vep$, and
the antipode $S$ defined by
\begin{align}
\begin{aligned}
& \Del(q^h) = q^h \ot q^h, \\
& \Del(e_i) = e_i \ot K_i^{-1} + 1 \ot e_i,
\quad \Del(f_i) = f_i \ot 1 + K_i \ot f_i, \\
& \vep(q^h)=1, \quad
\vep(e_i) = \vep(f_i) = 0, \\
& S(q^h) = q^{-h}, \quad
S(e_i) = - e_i K_i, \quad
S(f_i) = - K_i^{-1} f_i
\end{aligned}
\end{align}
for $h \in \pv$, $i \in I$ (see, for example, \cite{BKM98, Ka95}).

\smallskip

Let $U_q^+(\g)$ (resp.\ $\Um$) be the subalgebra of $\Uq$
generated by the elements $e_i$ (resp.\ $f_i$) for $i \in I$,
and let $U_q^0(\g)$ be the subalgebra of $\Uq$ generated by
$q^h$ $(h \in \pv)$.
Then we have the {\it triangular decomposition}
(\cite{BKM98, Ka95})\,:
\begin{align}
\Uq \cong \Um  \ot U_q^0(\g)  \ot U_q^+(\g).
\end{align}

\begin{rem}
In \cite{Ka95}, Kang considered the generalized Kac-Moody
algebras associated with Borcherds-Cartan matrices of
{\it charge} $\m = (m_i \in \Z_{>0} \mid i\in I, \ \ m_i=1
\ \mbox{ for } \ i \in I^\re)$.
The charge $m_i$ is the multiplicity of the simple root
corresponding to $i \in I$.
For example, the {\it Monster Lie algebra}, which played an
important role in Borcherds' proof of Monstrous moonshine,
is a generalized Kac-Moody algebra associated with the
Borcherds-Cartan matrix $A=(-(i+j))_{i,j\in I}$ of
charge $\m=(c(i) \mid i\in I)$, where
$I=\{ -1 \} \cup \{ 1, 2, \dotsc \}$ and $c(i)$ are the
coefficients of the elliptic modular function
\[ J(q) = j(q)-744 = \sum_{n=-1}^{\infty} c(n) q^n
= q^{-1} + 196884 q + 21493760 q^2 + \dotsb. \]
In this paper, we assume that $m_i = 1$ for all $i\in I$.
However, we do not lose generality by this assumption.
Indeed, if we take Borcherds-Cartan matrices with some of the rows and
columns identical, then the generalized Kac-Moody algebras
with charge introduced
in \cite{Ka95} can be recovered from the ones in this paper
by identifying the $h_i$'s and $d_i$'s (hence $\al_i$'s)
corresponding to these identical rows and columns.
\end{rem}

\medskip

\section{{\bf Representation Theory}}

A $\Uq$-module $V$ is called a {\it weight module}
if it admits a {\it weight space decomposition}
$V = \bop_{\mu \in P} V_\mu$, where
$V_\mu \seteq \set{v \in V}
{q^h v = q^{\mu(h)} v \mbox{ for every } h \in \pv}$.
We call $\wt(V) \seteq \set{\mu \in P}{V_\mu \ne 0 }$
the set of {\it weights} of $V$.

Let $\CO$ be the category consisting of
weight modules $V$
with finite-dimensional weight spaces such that
$\wt(V) \subset \bigcup_{j=1}^s (\la_j - Q_+)$
for finitely many $\la_1, \dotsc, \la_s \in P$.

One of the most interesting examples of $\Uq$-modules in the
category $\CO$ is the class of highest weight modules
defined below.
A weight module $V$ over $\Uq$ is called a
{\it highest weight module with highest weight $\la \in P$}
if there exists a {\em non-zero} vector $v \in V$
(called a {\it highest weight vector\/}) such that
\begin{tenumerate}
\item $V= \Uq v$,
\item $v\in V_\la$,
\item $e_i \, v = 0$ for every $i \in I$.
\end{tenumerate}

Let $J(\la)$ denote the left ideal of
$\Uq$ generated by $e_i$,
$q^h - q^{\la(h)}$ ($i \in I$, $\ h \in \pv$), and set
$M(\la) = \Uq / J(\la)$.
Then, via left multiplication,
$M(\la)$ becomes a $\Uq$-module
called the {\it Verma module}.
The basic properties of Verma modules are summarized
in the following proposition.

\begin{prop}[\cite{BKM98, Ka95}] \hfill
\begin{enumerate}
\item $M(\la)$ is a highest weight $\Uq$-module
with highest weight $\la$ and highest weight vector
$v_\la \seteq 1 + J(\la)$.
\item $M(\la)$ is a free $\Um$-module
generated by $v_\la$.
\item Every highest weight $\Uq$-module
with highest weight $\la$ is a quotient of $M(\la)$.
\item $M(\la)$ contains a unique maximal submodule $R(\la)$.
\end{enumerate}
\end{prop}

\noindent
The irreducible quotient $\vla = M(\la) / R(\la)$ is
an {\it irreducible highest weight $\Uq$-module
with highest weight $\la$}.

In the rest of this section, we focus on the structure
of the irreducible highest weight $\Uq$-module $\vla$
with a dominant integral highest weight $\la \in P^+$.

We first recall:

\begin{prop}[\cite{BKM98, Ka95}] \label{prop:hw-module}
Let $\la \in P^+$ be a dominant integral weight.
\begin{enumerate}
\item The highest weight vector
$v_\la$ of $\vla$ satisfies the following relations
\begin{align}\label{eq:highest}
\begin{cases}
f_i^{\la(h_i)+1} v_\la =0 \quad & \mbox{if }\ a_{ii}=2, \\
f_i v_\la =0 & \mbox{if }\ \la(h_i)=0.
\end{cases}
\end{align}
\item Conversely, let $V$ be a highest weight $\Uq$-module
with a highest weight $\la \in P^+$ and a highest weight vector $v$.
If $v$ satisfies the relations in \eqref{eq:highest},
then $V$ is isomorphic to $\vla$.
\end{enumerate}
\end{prop}

Consider the anti-involution of $\Uq$ defined by
\begin{align}
q^h \mt q^h, \quad
e_i \mt q_i f_i K_i^{-1}, \quad
f_i \mt q_i^{-1} K_i e_i
\quad (i \in I).
\end{align}
By standard arguments, we see that there exists
a unique symmetric bilinear form
$(\ , \ )$ on $\vla$ with $\la \in P^+$ satisfying
\begin{align} \label{eq:bilinear form}
\begin{aligned}
& (v_\la, v_\la) = 1, \quad (q^h u, v) = (u, q^h v), \\
& (e_i u, v) = (u, q_i f_i K_i^{-1} v),
\quad (f_i u, v) = (u, q_i^{-1} K_i e_i v)
\end{aligned}
\end{align}
for $i \in I$, $h \in \pv$, $u, v \in \vla$.

Note that $(\vla_\mu, \vla_{\tau})=0$ unless $\mu=\tau$.
Moreover, we have:

\begin{lem}
The symmetric bilinear form $(\ , \ )$ on $\vla$ defined by
\eqref{eq:bilinear form} is non-degenerate.
\end{lem}
\begin{proof}
The kernel $\set{u \in \vla}{(u, \vla) = 0}$
is a $\Uq$-submodule not containing $v_\la$.
Since $\vla$ is irreducible, it must vanish.
\end{proof}

The weights of the irreducible highest weight $\Uq$-module
$\vla$ satisfy the following properties.

\begin{prop} \label{prop:V(la)}
Let $\mu$ be a weight of $\vla$ with $\la \in P^+$,
and let $i \in I^\im$.
\begin{enumerate}
\item $\mu(h_i) \in \Z_{\ge 0}$.
\item If $\mu(h_i)=0$, then $\vla_{\mu-\al_i}=0$.
\item If $\mu(h_i)=0$, then $f_i (\vla_\mu) =0$.
\item
If $\mu(h_i)\le 2c_i$, then $e_i(\vla_\mu)=0$.

\end{enumerate}
\end{prop}
\begin{proof}
Set $\mu = \la-\al$ with $\al \in Q_+$.
Write $\al = \sum_{\nu=1}^r \al_{i_\nu}$.
Since $\al_{i_\nu}(h_i) \le 0$, we have $\al(h_i)\le 0$,
which yields $\mu(h_i)=\la(h_i)-\al(h_i) \ge \la(h_i)\ge 0$.
Note that $\mu(h_i)=0$ implies that $\la(h_i) = \al_{i_\nu}(h_i)=0$
$(\nu = 1, \dotsc, r)$.
Thus $f_i$ commutes with all the $f_{i_\nu}$'s,
and kills $v_\la$. Therefore we obtain $\vla_{\mu-\al_i}=0$.

The properties (c) and (d) easily follow from (a) and (b).
\end{proof}

Fix $i \in I$, and let
$U_i$ be the subalgebra of $\Uq$ generated by
$e_i$, $f_i$ and $U_q^0(\g)$.
Then the algebra $U_i$ can be viewed as a rank 1 quantum generalized
Kac-Moody algebra associated with
the Borcherds-Cartan matrix $A=(a_{ii})$.
For example, if $a_{ii}=0$,
then $U_i$ is isomorphic to the {\it quantum Heisenberg algebra}.

Let $V= \vla$ be the irreducible highest weight $U_i$-module
with a dominant integral highest weight $\la$.
If $a_{ii}=2$, it is well-known that
$\dim V = \la(h_i) + 1$.
If $a_{ii}\le 0$, then $\dim V$ is given in the following lemma, which
easily follows
from Proposition \ref{prop:hw-module}.

\begin{lem} \label{lem:rep of U_i}
Suppose $a_{ii} \le 0$ and
let $\vla$ be the irreducible highest weight $U_i$-module
with a dominant integral highest weight $\la$.
Then we have
$\vla=\Q(q)v_\la$ if $\la(h_i)=0$, and
$\{ f_i^n v_\la \}_{n \ge 0}$ is a basis of $\vla$ if $\la(h_i) >0$.
\end{lem}

Observe that the $U_i$-module structure on $\vla$ is given by
\begin{align} \label{eq:Ui action}
\begin{aligned}
& K_if_i^n v_\la = q_i^{\la(h_i)+2nc_i}
f_i^n v_\la, \\
& f_i f_i^n v = f_i^{n+1} v, \\
& e_i f_i^n v_\la =
\{ n \}_i [\la(h_i)+c_i (n-1)]_i f_i^{n-1} v_\la.
\end{aligned}
\end{align}

\medskip

\section{{\bf Category $\CO_{int}$}}

\begin{defi} \label{defi:O_{int}}
The {\it category $\CO_{int}$} consists of $\Uq$-modules
$M$ satisfying the following properties\,:
\begin{tenumerate}
\item $M=\bop_{\mu \in P} M_\mu$ with
$\dim M_\mu < \infty$ for every $\mu \in P$, where
\begin{align*}
M_\mu \seteq \set{v\in M}
{q^h v = q^{\mu(h)} v \mbox{ for every }h \in \pv},
\end{align*}
\item there exist a finite number of weights
$\la_1, \dotsc, \la_s \in P$ such that
\[ \wt(M) \seteq \set{\mu \in P}{M_\mu \ne 0}
\subset \bigcup_{j=1}^s \left(\la_j - Q_+\right), \]
\item if $a_{ii}=2$, then the action of $f_i$ on $M$ is locally nilpotent,
\item if $a_{ii}\le 0$, then $\mu(h_i) \in \Z_{\ge 0}$
for every $\mu \in \wt(M)$,
\item if $a_{ii} \le 0$ and $\mu(h_i) = 0$, then
$f_i M_\mu = 0$.
\item
if $a_{ii} \le 0$ and $\mu(h_i) =2c_i$, then
$e_i M_\mu = 0$.
\end{tenumerate}
\end{defi}

\begin{rem} \hfill
\begin{tenumerate}
\item
Note that if $a_{ii}\le 0$ and $\mu(h_i)\le 2c_i$ then
$e_iM_\mu=0$. Indeed, if $\mu(h_i)<2c_i$, then $(\mu+\al_i)(h_i)<0$ and
$\mu+\al_i$ is not a weight of $M$ by Definition \ref{defi:O_{int}} (iv).
\item
By (ii), the action of $e_i$ on a $\Uq$-module in the category $\CO_{int}$
is locally nilpotent.
\item
Note that a submodule or a quotient module of a
$\Uq$-module in the category $\CO_{int}$ also
belongs to the category $\CO_{int}$.
Furthermore, it is straightforward to verify that the category
$\CO_{int}$ is closed under taking
a finite number of direct sums and tensor products.
\end{tenumerate}
\end{rem}

\begin{prop} \label{prop:P^+}
Let $\vla$ be the irreducible highest weight $\Uq$-module
with highest weight $\la \in P$.
Then $\vla$ belongs to the category $\CO_{int}$
if and only if $\la \in P^+$.
\end{prop}
\begin{proof}
As we have seen in Proposition \ref{prop:V(la)},
$\vla$ belongs to the category $\CO_{int}$
for $\la \in P^+$.

Conversely, assume that $\vla$ lies in the category $\CO_{int}$.
If $a_{ii}=2$, since $f_i$ is locally nilpotent on $\vla$,
there exists a non-negative integer $N_i$ such that
$f_i^{N_i} v_\la \ne 0$ and $f_i^{N_i + 1} v_\la =0$.
Hence we have
\begin{align*}
0= e_i f_i^{N_i + 1} v_\la
= [ N_i + 1 ]_i [\la(h_i) - N_i]_i f_i^{N_i} v_\la,
\end{align*}
which implies $\la(h_i) = N_i \in \Z_{\ge 0}$.

If $a_{ii}\le 0$,
we have $\la(h_i) \in \Z_{\ge 0}$ by Definition \ref{defi:O_{int}} (iv).
\end{proof}

\smallskip

\begin{prop} \label{prop:V=V(la)} \hfill
\begin{enumerate}
\item If $V$ is a highest weight $\Uq$-module in the category
$\CO_{int}$ with highest weight $\la \in P$,
then $\la \in P^+$ and $V \simeq \vla$.
\item Every irreducible $\Uq$-module in the category
$\CO_{int}$ is isomorphic to some $\vla$
with $\la \in P^+$.
\end{enumerate}
\end{prop}
\begin{proof}
(a)\quad Let $v$ be a highest weight vector of $V$ with weight $\la$.
If $a_{ii}=2$, by the proof of Proposition \ref{prop:P^+},
we have $\la(h_i) \in \Z_{\ge 0}$
and $f_i^{\la(h_i)+1} v=0$.

If $a_{ii}\le 0$, then by the definition of $\CO_{int}$,
we have $\la(h_i) \in \Z_{\ge 0}$ and $\la(h_i)=0$ implies
$f_i v=0$.
Hence $\la \in P^+$, and
Proposition \ref{prop:hw-module} (b) implies that
$V \simeq \vla$.

\smallskip
\noindent
(b)\quad Let $V$ be an irreducible $\Uq$-module in the category
$\CO_{int}$.
Since $\wt(V) \subset \bigcup_{j=1}^s (\la_j - Q_+)$ for some
$\la_j \in P$ $(j=1, \dotsc, s)$, there exists a
non-zero {\it maximal vector} $v\in V_\la$
for some $\la \in \wt(V) \subset P$
such that $e_i v= 0$ for every $i \in I$.
Then $V=\Uq v$ is isomorphic to $\vla$ by (a).
\end{proof}

We will now prove that every $\Uq$-module in the category
$\CO_{int}$ is semisimple.
Let $\vp$ be the anti-involution of $\Uq$ given by
\begin{align}
e_i\mt f_i,\quad f_i\mt e_i,\quad q^h\mt q^h.
\end{align}

Let $M=\bop_{\mu \in P} M_\mu$ be a $\Uq$-module
in the category $\CO_{int}$.
We define its {\it finite dual} to be the vector space
\begin{align}
M^* \seteq \bop_{\mu \in P} M_\mu^*, \mbox{ where }
M_\mu^* \seteq \operatorname{Hom}_{\Q(q)} (M_\mu, \Q(q)).
\end{align}
Using the anti-involution $\vp$ of $\Uq$, we define a $\Uq$-module
structure on $M^*$ by
\begin{align}
\lan x \cdot \phi,\ v \ran
= \lan \phi,\ \vp(x) \cdot v \ran
\quad \mbox{for} \ \ x \in \Uq, \ \phi \in M^*,
v \in M.
\end{align}

The following lemma is an immediate consequence of the definitions.

\begin{lem} \label{lem:dual O_{int}}
Let $M=\bop_{\mu \in P} M_\mu$ be a $\Uq$-module
in the category $\CO_{int}$.
\begin{enumerate}
\item There exists a canonical isomorphism
$(M^*)^* \cong M $ as $\Uq$-modules.
\item The space $M_\mu^*$ is the weight space of $M^*$ with weight $\mu$.
\item We have
\begin{align*}
\wt(M^*) =\wt(M).
\end{align*}
\item $M^*$ also belongs to the category $\CO_{int}$.
\end{enumerate}
\end{lem}
\begin{proof}
Since the other statements are evident,
we only show that $M^*$ satisfies (v) and (vi)
in Definition \ref{defi:O_{int}}.

\smallskip
\noindent
(v)\quad For $i \in I^\im$ and $\mu\in P$ with $\mu(h_i)=0$, we have
\[ \lan f_i M^*_\mu, M \ran = \lan M^*_\mu, e_i M_{\mu-\al_i}\ran=0, \]
where the last equality follows from $\lan h_i,\mu-\al_i\ran=2c_i$.

\smallskip
\noindent
(vi)\quad For $i \in I^\im$ and $\mu\in P$ with $\mu(h_i)=2c_i$, we have
\[ \lan e_iM^*_\mu, M\ran=\lan M^*_\mu, f_iM_{\mu+\al_i}\ran=0, \]
where the last equality follows from $\lan h_i,\mu+\al_i\ran=0$.
\end{proof}

\smallskip

Let $M$ be a non-zero $\Uq$-module in the category $\CO_{int}$
and choose a {\it maximal weight} $\la \in \wt(M)$ with the
property that $\la+ \al_i$ is {\it not} a weight of $M$
for any $i \in I$.
Fix a non-zero vector $v_\la \in M_\la$ and set
$V = \Uq v_\la$.
Then by Proposition \ref{prop:V=V(la)}, we have $\la \in P^+$
and $V$ is isomorphic to $\vla$.

Let us take $\phi_\la\in M^*_\la$
satisfying $\phi_\la(v_\la)=1$, and set
$W = \Uq \phi_\la$.
Then, again by Proposition \ref{prop:V=V(la)},
$W$ is isomorphic to $\vla$.

\begin{lem} \label{lem:one factor}
Let $M$ be a $\Uq$-module in the category $\CO_{int}$
and let $V$ be the submodule of $M$ generated by a maximal vector
$v_\la$ of weight $\la$.
Then we have
\begin{align*}
M \simeq V \op (M/V).
\end{align*}
\end{lem}
\begin{proof}
We will show that the short exact sequence
\begin{align*}
0 \longrightarrow V \longrightarrow M
\longrightarrow M/V \longrightarrow 0
\end{align*}
splits.

Take the dual of the inclusion $W \hookrightarrow M^*$
to get a map $(M^*)^* \longrightarrow W^*$.
Since $M \cong (M^*)^*$, we obtain a map
$\eta\cl M \isoto(M^*)^* \longrightarrow W^*$,
which yields
\begin{align*}
\psi \cl V \hookrightarrow M\xrightarrow{\ \eta \ }W^*.
\end{align*}
It is easy to check that the image of $v_\la$ under $\psi$ is non-zero.
Since $V \simeq \vla$ and
$W^* \simeq \vla^* \simeq \vla$, $\psi$ is an isomorphism
by Schur's Lemma.
Hence the above short exact sequence splits,
and we have $M \simeq V \op (M/V)$.
\end{proof}

Using this lemma, we will prove the complete reducibility theorem
for the category $\CO_{int}$.

\begin{thm} \label{thm:complete reducibility}
Every $\Uq$-module in the category $\CO_{int}$ is isomorphic
to a direct sum of irreducible highest weight modules $\vla$ with
$\la \in P^+$.
\end{thm}
\begin{proof}
By Proposition \ref{prop:V=V(la)} (b),
it is enough to show that any $\Uq$-module $M$
in the category $\CO_{int}$ is semisimple.

Let $U^{\ge 0}$ denote the subalgebra of $\Uq$
generated by the elements $q^h$ $(h \in \pv)$ and
$e_i$ $(i \in I)$.
We shall first show that, if $M$ is generated
as a $\Uq$-module by a finite-dimensional
$U^{\ge 0}$-submodule $F$, then
$M$ is semisimple.
We shall argue by the induction on the dimension of $F$.
If $F \ne 0$, let us choose a maximal weight vector $v_\la$
of weight $\la \in P^+$ in $F$, and set $V=\Uq v_\la\simeq \vla$.
By Lemma \ref{lem:one factor}, we have
\begin{align*}
M\simeq V \op (M/V).
\end{align*}
Since $M/V = \Uq \bigl(F/(F\cap V)\bigr)$
and $\dim \bigl(F/(F\cap V)\bigr) < \dim F$,
the induction hypothesis implies that $M/V$
is semisimple.
Hence $M$ is semisimple.

Now let $M$ be an arbitrary $\Uq$-module in the category $\CO_{int}$.
Then, for any $v\in M$, $U^{\ge 0}v$ is finite-dimensional, and hence
$\Uq v$ is semisimple.
Since $M$ is a sum of semisimple $\Uq$-modules $\Uq v$,
a standard argument for semisimplicity
(\cite[Proposition 3.12]{CR81}) implies that
$M$ is semisimple.
\end{proof}

\smallskip

\begin{coro} \label{coro:i-string}
Let $M$ be a $\Uq$-module in the category $\CO_{int}$.
Then for any $i \in I$, $M$ satisfies the following properties:
\begin{enumerate}
\item If $a_{ii}=2$, then $M$ is
a direct sum of finite-dimensional irreducible $U_i$-submodules.
\item If $a_{ii}\le 0$, then $M$ is a direct sum of 1-dimensional
or infinite-dimensional irreducible highest weight
$U_i$-submodules.
\end{enumerate}
\end{coro}
\begin{proof}
Our assertions immediately follow
from Lemma \ref{lem:rep of U_i} and
Theorem \ref{thm:complete reducibility}.
\end{proof}

\medskip

\section{{\bf Crystal Bases}}

In this section, we will develop the {\it crystal basis theory} for
$\Uq$-modules in the category $\CO_{int}$.
The following fundamental lemma easily follows from
Corollary \ref{coro:i-string}.
\begin{lem}\label{lem:str}
Let $M$ be a $\Uq$-module
in the category $\CO_{int}$.
Then for any $i \in I$, we have\,{\rm :}
\begin{enumerate}
\item $M=\bop_{n\in\Z_{\ge 0}}f_i^{(n)}(\ker e_i)$.
\item For any $\la\in P$ and $n > 0$, the map
\begin{align}\label{eq:fn}
f_i^{(n)}\cl \ker e_i\cap M_\la\to M_{\la-n\al_i}
\end{align}
is either a monomorphism or zero.
More precisely, the morphism \eqref{eq:fn} is injective
if $i \in I^\re$ and $\la(h_i)\ge n$ or if $i \in I^\im$
and $\la(h_i)>0$.
Otherwise, it vanishes.
\end{enumerate}
\end{lem}

Hence for any weight vector $u \in M_\la$
there exists a unique family $\{ u_n \}_{n \in \Z_{\ge 0}}$
of elements of $M$ such that
\begin{align}\label{eq:str}
\parbox{300pt}{
\begin{enumerate}
\item
$u_n\in\ker e_i\cap M_{\la+n\al_i}$,
\vspace{5pt}
\item
$u = \sum_{n\ge 0}f_i^{(n)}u_n$,
\vspace{5pt}
\item
if $i \in I^\re$ and $\lan h_i,\la+n\al_i\ran < n$,
then $u_n=0$,
\vspace{5pt}
\item
if $i \in I^\im$, $n>0$ and $\lan h_i,\la+n\al_i\ran=0$,
then $u_n=0$.
\end{enumerate}
}
\end{align}
Note that $u_n=0$ for $n \gg 0$, because
$\la+n\al_i$ is not a weight of $M$ for $n \gg 0$.
We call the expression (b) in \eqref{eq:str}
the {\it \str} of $u$.

\begin{rem}\hfill
\begin{tenumerate}
\item
The conditions (c) and (d) in \eqref{eq:str}
are equivalent to saying that
$f_i^{(n)}u_n=0$ implies $u_n=0$.
\item
If $\sum_{n\ge 0}f_i^{(n)}u_n$ is an \str,
then
$\sum_{n\ge 0}f_i^{(n)}u_{n+1}$ is an \str.
However,
$\sum_{n\ge 1}f_i^{(n)}u_{n-1}$ is not necessarily an \str.
\end{tenumerate}
\end{rem}

\begin{defi}
For $i \in I$ and $u\in M$,
let $u =\sum_{n\ge 0} f_i^{(n)} u_{n}$
be the \str\ of $u$.
Then the {\it Kashiwara operators}
$\eit$ and $\fit$ on $M$ are defined by
\begin{align}
\eit \, u = \sum_{n\ge 1} f_i^{(n-1)} u_n, \qquad
\fit \, u = \sum_{n\ge 0} f_i^{(n+1)} u_n.
\end{align}
\end{defi}

By Lemma \ref{lem:str},
if $u =\sum_{n\ge 0} f_i^{(n)} u_n$ with $u\in M_\la$ and
$u_n \in \ker e_i \cap M_{\la+n\al_i}$, then we have
\begin{align}\label{eq:te}
\eit \, u = \sum_{\substack{n\ge 1,\\ f_i^{(n)}u_n\ne 0}}
 f_i^{(n-1)} u_n \quad\mbox{and}\quad
\fit \, u = \sum_{n \ge 0} f_i^{(n+1)} u_n.
\end{align}

Note that if $a_{ii}\le 0$, then $\fit = f_i$,
but $\eit \ne e_i$ (see Lemma \ref{lem:EF} below and the
remark after \eqref{eq:qKe}).
Note that we have
\begin{align}\label{eq:ef}
\eit\circ\fit \mid_{M_\la}=\id_{M_\la}
\mbox{ if } a_{ii}\le 0 \mbox{ and }\la(h_i)>0.
\end{align}

\smallskip

We now define the notion of {\it crystal bases} for $\Uq$-modules
in the category $\CO_{int}$.
Let $\Az$ be the localization of the polynomial ring $\Q[q]$
at the prime ideal $\Q[q]q$:
\begin{align*}
\Az = \set{f/g}{f, g \in \Q[q], \ g(0) \ne 0} .
\end{align*}
Thus the local ring $\Az$ is a principal ideal domain with
$\Q(q)$ as its field of quotients.

\smallskip

\begin{defi}\label{def:cl}
Let $M$ be a $\Uq$-module in the category $\CO_{int}$.
A {\it crystal lattice} of $M$ is a free $\Az$-submodule $L$ of $M$
satisfying the following conditions\,:
\begin{tenumerate}
\item $L$ generates $M$ as a vector space over $\Q(q)$,
\item $L=\bop_{\la \in P} L_\la$, where
$L_\la \seteq L \cap M_\la$,
\item $\eit \, L \subset L$ and
$\fit \, L \subset L$ \ for every $i \in I$.
\end{tenumerate}
\end{defi}

The map $\Az\ni f \longmapsto f(0)\in\Q$ induces an isomorphism
\begin{align*}
\Az / q\Az\isoto \Q,
\end{align*}
and we have
\begin{align*}
\Q \ot_{\Az} L \isoto L/qL.
\end{align*}
The passage from $L$ to the quotient $L/qL$ is referred to as
taking the {\it crystal limit}.
Since the Kashiwara operators $\eit$ and $\fit$ preserve the
crystal lattice $L$, they also define operators on $L/qL$,
which we will denote by the same symbols.

\smallskip

\begin{lem}\label{lem:EF}
Let $L$ be a free $\Az$-submodule $L$ of $M\in\CO_{int}$
satisfying the conditions {\rm(i)} and {\rm(ii)}
in Definition $\ref{def:cl}$.
\begin{enumerate}
\item $L$ is a crystal lattice if and only if
$L=\bop_{n \ge 0}f_i^{(n)}(\ker e_i\cap L)$ for every $i \in I$.
\item Let $E_i$ and $F_i$ be operators on $M$ such that
\[ E_i f_i^{(n)}v=a_{i,n,\la} f_i^{(n-1)}v, \quad
F_i f_i^{(n)}v = b_{i,n,\la} f_i^{(n+1)}v \]
for $v\in \ker e_i\cap M_\la$ with $f_i^{(n)}v \ne 0$.
Here, $a_{i,n,\la}$ and $b_{i,n,\la}$ are invertible elements of $\Az$.

Then $L$ is a crystal lattice if and only if
$E_i L \subset L$ and $F_i L \subset L$.
Moreover, if we define the operators
$R_i$ and $S_i$ by
\[ R_i f_i^{(n)}v = a_{i,n,\la}(0)f_i^{(n)}v, \quad
S_i f_i^{(n)}v = b_{i,n,\la}(0) f_i^{(n)}v, \]
then $R_i L \subset L$, $S_i L \subset L$ and
the induced action of
$E_i$ {\rm(}resp.\ $F_i${\rm)} on $L/qL$ coincides with
that of $\eit\circ R_i$ {\rm(}resp.\ $\fit\circ S_i${\rm)}.
In particular, we have
$E_{i}^{-1} (L)=\eit^{-1}(L)$.
\end{enumerate}
\end{lem}
\begin{proof}
Since the other statements are obvious, let us show that
if $E_i L \subset L$ and $F_i L \subset L$, then
$L=\bop_{n\ge 0}f_i^{(n)}(\ker e_i\cap L)$.
The condition $F_i L\subset L$ implies
that $L\supset \bop_{n\ge 0}f_i^{(n)}(\ker e_i\cap L)$.
In order to show the inverse inclusion,
let $\sum_{n=0}^Nf_i^{(n)}u_n$ be the \str\ of $u\in L_\la$.
Let us show that every $u_n$ belongs to $L$ by the induction of $N$.
Since we have
\[ E_iu=\sum_{n=1}^{N}a_{i,n,\la+n\al_i}f_i^{(n-1)}u_n\in L, \]
the induction hypothesis implies
$u_n\in L$ for $n\ge 1$.
Hence $u_0=u-\sum_{n\ge 1}f_i^{(n)}u_n$ belongs to $L$.
\end{proof}

Note that for $u\in \ker e_i\cap M_\la$, we have
\begin{align}\label{eq:qKe}
q_i^{-1}K_i e_i(f_i^{(n)}u)&=
\begin{cases}
q_i^{1-n}\dfrac{1-q_i^{2(a+1-n)}}{1-q_i^2}f_i^{(n-1)}u
&\mbox{if $i \in I^\re$,}\\
\dfrac{(1-q_i^{2nc_i})(1-q_i^{2(a+(n-1)c_i)})}{(1-q_i^2)(1-q_i^{2c_i})}
f_i^{(n-1)}u
&\mbox{if $a_{ii}<0$,}\\
n\dfrac{1-q_i^{2a}}{1-q_i^2}f_i^{(n-1)}u&\mbox{if $a_{ii}=0$,}
\end{cases}
\end{align}
where $a=\la(h_i)$.
Hence when $a_{ii}\le 0$, we can use $q_i^{-1}K_i e_i$ instead of $\eit$.

\begin{defi}
Let $M$ be a $\Uq$-module in the category $\CO_{int}$.
A {\it crystal basis} of $M$ is a pair $(L, B)$ such that
\begin{tenumerate}
\item $L$ is a crystal lattice of $M$,
\item $B$ is a $\Q$-basis of
$L/qL \cong \Q \ot_{\Az} L$,
\item $B = \bigsqcup_{\la \in P} B_\la$, where
$B_\la \seteq B \cap (L_\la/ qL_\la)$,
\item $\eit B \subset B \cup \ze$ and
$\fit B \subset B \cup \ze$ for every $i \in I$,
\item for any $i \in I$ and $b, b'\in B$, we have
$\fit b = b'$ if and only if $b=\eit b'$.
\end{tenumerate}
\end{defi}

The set $B$ is endowed with a colored oriented graph structure whose
arrows are defined by
\begin{align*}
b \stackrel {i} \longrightarrow b' \quad \mbox{if and only if}
\quad \fit b = b'.
\end{align*}
The graph $B$ is called the {\it crystal graph} of $M$
and it reflects the combinatorial structure of $M$.
For instance, we have
\begin{align*}
\dim_{\Q(q)} M_{\la}
=\mbox{rank}_{\Az} L_\la
=\# B_{\la}
\quad \mbox{for every }\ \la \in P.
\end{align*}
(See, for example, \cite[Theorem 4.2.5]{HK2002}.
See also Proposition \ref{prop:fn} below.)

\begin{thm} \label{thm:l_0}
Let $M$ be a $\Uq$-module
in the category $\CO_{int}$
and let $(L, B)$ be a crystal basis of $M$.
For $i \in I$ and $u\in L_\la$,
let $u = \sum_{l\ge 0} f_i^{(l)} u_l$ be the \str\ of $u$.
Then the following statements are true.
\begin{enumerate}
\item $u_l \in L$ for every $l\in \Z_{\ge 0}$.
\item If $u+qL \in B$, then there exists a non-negative integer
$l_0$ such that
\begin{enumerate}
\item $u_l \in q L$ for every $l \ne l_0$,
\item $u_{l_0} + qL \in B$,
\item $u \equiv f_i^{(l_0)} u_{l_0} \mod q L$.
\end{enumerate}
\end{enumerate}
In particular, we have
$\eit u \equiv f_i^{(l_0 - 1)} u_{l_0}$ and
$\fit u \equiv f_i^{(l_0 + 1)} u_{l_0} \mod q L$.
\end{thm}
\begin{proof}
The proof is the same as the one given
in \cite[Proposition 2.3.2]{Kas91} (see also
\cite[Proposition 4.2.11]{HK2002}).
\end{proof}

Let $M$ be a $\Uq$-module
in the category $\CO_{int}$ with a crystal basis $(L, B)$.
For $i \in I$, we define
the maps $\wt\cl B\to P$,
$\vep_i\cl B\to\Z_{\ge 0}$
and $\vp_i\cl B \to \Z_{\ge 0} \cup \{ \infty \}$ by
\begin{align}
\begin{aligned}
& \wt(b) = \la \quad \mbox{ for } b \in B_\la \ (\la \in P),\\
& \vep_i (b) = \max \set{l \ge 0}{\eit{}^l b \in B}, \\
& \vp_i (b) = \max \set{l \ge 0}{\fit{}^l b \in B}.
\end{aligned}
\end{align}
Note that for $i\in I^\im$, we have (see Example \ref{exp:mod})
\begin{align*}
\vp_i(b)=\begin{cases}
0&\mbox{if } \lan h_i,\wt(b)\ran=0,\\
\infty&\mbox{if } \lan h_i,\wt(b)\ran>0.
\end{cases}
\end{align*}

Then, using Theorem \ref{thm:l_0}, we can easily prove the
following proposition.

\begin{prop} \hfill
\begin{enumerate}
\item If $b\in B_\la$ satisfies $\eit b \in B$, then
\begin{align*}
\wt(\eit b) = \la + \al_i, \quad
\vep_i(\eit b) = \vep_i(b) - 1, \quad
\vp_i(\eit b) = \vp_i(b) + 1.
\end{align*}
\item If $b\in B_\la$ satisfies $\fit b \in B$, then
\begin{align*}
\wt(\fit b) = \la - \al_i, \quad
\vep_i(\fit b) = \vep_i(b) + 1, \quad
\vp_i(\fit b) = \vp_i(b) - 1.
\end{align*}
\end{enumerate}
\end{prop}

\begin{prop}\label{prop:fn}
For any $n \in \Z_{\ge 0}$ and $i \in I$, we have
\[ (f_i^n M \cap L)/ (f_i^n M \cap qL)
=\bop_{\substack{b\in B\\ \vep_i(b)\ge n}}\Q b. \]
In particular, we have
\[\dim_{\Q(q)} (f_i^n M)_\la= \# \set{b\in B_\la}{\vep_i(b)\ge n}.\]
\end{prop}
\begin{proof}
Since $L = \bop_{m\ge 0} f_i^{(m)} (\ker e_i \cap L)$
by Lemma \ref{lem:EF}, one has
\[ L/qL = \bop_{m\ge 0}
f_i^{(m)} (\ker e_i \cap L) / q f_i^{(m)} (\ker e_i \cap L) . \]
If $b\in B$ satisfies $\vep_i (b) = m$,
then $b$ belongs to
$f_i^{(m)} (\ker e_i \cap L) / q f_i^{(m)} (\ker e_i \cap L)$.
Hence $\set{b \in B}{\vep_i (b) = m}$
is a basis of
$f_i^{(m)} (\ker e_i \cap L) / q f_i^{(m)} (\ker e_i \cap L)$.
To complete the proof, it is enough to note that
$f_i^n M \cap L = \op_{m \ge n} f_i^{(m)} (\ker e_i \cap L)$.
\end{proof}

\medskip

\section{{\bf Tensor Product Rule}}

Fix $i \in I$, and let
$U_i$ be the subalgebra of $\Uq$ generated by
$e_i$, $f_i$ and $U_q^0(\g)$.
In the following example, we investigate the structure of
crystal bases for irreducible highest weight $U_i$-modules.

\begin{ex}\label{exp:mod}
Let $V \seteq \vla$ be the irreducible highest weight
$U_i$-module with $m \seteq \la(h_i)\in \Z_{\ge 0}$
and highest weight vector $v_\la$.

\smallskip

\begin{enumerate}
\item If $a_{ii}=2$, then we have
\begin{align*}
K_i v_\la = q_i^{m} v_\la, \quad
e_i v_\la = 0, \quad
f_i^{m+1} v_\la =0,
\end{align*}
and $V = \bop_{l=0}^m \Q(q) f_i^{(l)} v_\la$.
Set
\begin{align*}
L=\bop_{l=0}^m \Az f_i^{(l)} v_\la
\quad \mbox{and} \quad
B = \{ v_\la, f_i v_\la, \dotsc, f_i^{(m)} v_\la \}.
\end{align*}
Then $(L, B)$ is a crystal basis of $V$ with
the Kashiwara operators given by
\begin{align*}
\eit(f_i^{(l)} v_\la) = f_i^{(l-1)} v_\la, \quad
\fit(f_i^{(l)} v_\la) = f_i^{(l+1)} v_\la.
\end{align*}
Moreover, we have
\begin{align*}
\wt(f_i^{(l)} v_\la) = \la - l \al_i, \quad
\vep_i(f_i^{(l)} v_\la) = l, \quad
\vp_i(f_i^{(l)} v_\la) = m-l.
\end{align*}

\smallskip

\item If $a_{ii} \le 0$ and $\la(h_i)=0$, then we have
\begin{align*}
K_i v_\la = v_\la, \quad e_i v_\la = f_i v_\la =0,
\quad \mbox{and} \quad V=\Q(q) v_\la.
\end{align*}
Set $L=\Az v_\la$ and $B = \{ v_\la \}$.
Then $(L,B)$ is a crystal basis of $V$, and
\begin{align*}
\eit v_\la = \fit v_\la = 0, \quad
\wt(v_\la) = \la, \quad
\vep_i(v_\la)=\vp_i(v_\la) =0.
\end{align*}

\smallskip

\item If $a_{ii} \le 0$ and $m \seteq \la(h_i)>0$, then we have
\begin{align*}
K_i v_\la = q_i^m v_\la, \quad e_i v_\la =0, \quad
\mbox{and} \quad
V =\bop_{l\ge 0} \Q(q) f_i^l v_\la.
\end{align*}
Set
\begin{align*}
L=\bop_{l\ge 0} \Az f_i^l v_\la
\quad \mbox{and} \quad
B= \set{f_i^l v_\la}{l\ge 0}.
\end{align*}
Then $(L, B)$ is a crystal basis of $V$ with the Kashiwara operators
given by
\begin{align*}
\eit(f_i^l v_\la) = f_i^{l-1} v_\la,
\quad
\fit(f_i^l v_\la) = f_i^{l+1} v_\la.
\end{align*}
Moreover, we have
\begin{align*}
\begin{aligned}
\wt(f_i^l v_\la) = \la - l \al_i, \quad
\vep_i(f_i^l v_\la) = l, \quad
\vp_i(f_i^l v_\la) = \infty.
\end{aligned}
\end{align*}
\end{enumerate}
\end{ex}

Now we are ready to state the {\it tensor product rule} for
crystal bases of $\Uq$-modules in the category $\CO_{int}$.

\begin{thm} \label{thm:tensor product rule}
Let $M_j$ be a $\Uq$-module in the category $\CO_{int}$ and
let $(L_j, B_j)$ be a crystal basis of $M_j$ $(j=1, 2)$.
Set
\begin{align*}
M= M_1 \ot_{\Q(q)} M_2,
\quad
L = L_1  \ot_{\Az} L_2, \quad B = B_1 \ot B_2.
\end{align*}
Then $(L,B)$ is a crystal basis of $M$, where the
Kashiwara operators $\eit$ and $\fit$ $(i \in I)$
on $B$ are given as follows\,{\rm :}
\begin{align*}
\begin{aligned}
& \eit(b_1\ot b_2)
= \begin{cases}
\eit b_1 \ot b_2 \quad & \mbox{if } \ \vp_i(b_1) \ge \vep_i(b_2),\\
b_1 \ot \eit b_2 \quad & \mbox{if } \ \vp_i(b_1) < \vep_i(b_2),
\end{cases}\\
& \fit(b_1\ot b_2)
= \begin{cases}
\fit b_1 \ot b_2 \quad & \mbox{if } \ \vp_i(b_1) > \vep_i(b_2),\\
b_1 \ot \fit b_2 \quad & \mbox{if } \ \vp_i(b_1) \le \vep_i(b_2).
\end{cases}
\end{aligned}
\end{align*}
\end{thm}
\begin{proof}
By Theorem \ref{thm:l_0}, it suffices to prove
the tensor product rule for irreducible highest weight
$U_i$-modules $M_1=\vla$
and $M_2=\vmu$ with $a \seteq \la(h_i) \in \Z_{\ge 0}$
and $b \seteq \mu(h_i) \in \Z_{\ge 0}$.

If $a_{ii}=2$, our assertion was already proved in \cite{Kas90}.
(See also \cite[Theorem 4.4.3]{HK2002}, \cite[Th\'eor\`eme 2.3.5]{SMF}).

If $a_{ii} \le 0$, set $c_i=-\dfrac{1}{2} a_{ii} \in
\Z_{\ge 0}$.
If $a=0$ or $b=0$, then our assertion is trivial.
We assume that $a$, $b>0$. Then, by \eqref{eq:ef}, we have
$\eit\circ\fit=\id_{M_1 \ot M_2}$ and
\begin{align*}
\begin{array}{ll}
M_1 = \vla = \bop_{l\ge 0} \Q(q) f_i^l u,
&M_2 = \vmu = \bop_{l\ge 0} \Q(q) f_i^l v,\\[3pt]
L_1=\bop_{l\ge 0} \Az f_i^l u,
&L_2=\bop_{l\ge 0} \Az f_i^l v,\\
B_1=\set{f_i^l u \mod q L_1}{l \in \Z_{\ge 0}},
&B_2=\set{f_i^l v \mod q L_2}{l \in \Z_{\ge 0}},
\end{array}
\end{align*}
where
\begin{align*}
e_i u = e_i v = 0, \quad
K_i u = q_i^{\la(h_i)}u = q_i^{a} u,
\quad
K_i v = q_i^{\mu(h_i)}v = q_i^{b} v.
\end{align*}
Consider the $\Az$-lattice
\begin{align*}
L = L_1 \ot_{\Az} L_2 = \bop_{l, m \ge 0} \Az (f_i^l u \ot f_i^m v).
\end{align*}
We have
\begin{align*}
f_i (f_i^l u \ot f_i^m v)
&= (f_i \ot 1 + K_i \ot f_i)(f_i^l u \ot f_i^m v) \\
&= f_i^{l+1} u \ot f_i^m v + q_i^{a + 2 l c_i} f_i^l u \ot f_i^{m+1} v
\\
&\equiv f_i^{l+1} u \ot f_i^m v \mod q L.
\end{align*}
Therefore, we have
$\fit L\subset L$ and
\[\fit(f_i^l u \ot f_i^m v)
\equiv \fit(f_i^{l} u )\ot f_i^m v \mod q L. \]

On the other hand, we have (see Lemma \ref{lem:EF} and \eqref{eq:qKe})
\begin{align*}
q_i^{-1}K_i e_i(f_i^l u \ot f_i^m v)
&= (q_i^{-1}K_i e_i \ot 1 + K_i \ot q_i^{-1}K_i e_i)(f_i^l u \ot f_i^m v)
\\
&= (q_i^{-1}K_i e_i f_i^{l} u) \ot f_i^m v +
q_i^{a+2c_il}f_i^l u \ot q_i^{-1}K_i e_i f_i^{m} v\\
&\equiv
(q_i^{-1}K_i e_i f_i^{l} u) \ot f_i^m v \mod q L.
\end{align*}
Therefore, we have
$\eit L\subset L$, and
$\eit(u \ot f_i^m v)\equiv 0 \mod q L$.
For $l > 0$, we have
\[ \eit(f_i^l\ot f_i^mv) \equiv \eit\fit(f_i^{l-1}\ot f_i^mv)
= f_i^{l-1}\ot f_i^mv\mod qL. \]
Hence $(L_1 \ot L_2, B_1\ot B_2)$ is a crystal basis
with the desired action of $\eit$ and $\fit$.
\end{proof}

\medskip

\section{{\bf Crystal basis for $\Um$}}

In this section, we will define the notion of crystal basis for $\Um$
following \cite{Kas91}.
Although \cite{Kas91} treated the Kac-Moody case,
similar arguments can be applied to
the generalized Kac-Moody algebra case with
slight modifications, and we omit the details.

Fix $i \in I$. For any $P \in \Um$, there exist unique
$Q$, $R \in \Um$ such that
\begin{align}\label{eq:boson decomp}
e_i P - P e_i = \frac {K_i Q - K_i^{-1} R} {q_i - q_i^{-1}}.
\end{align}
We define the endomorphisms $e_i'$, $e_i''\cl\Um \to \Um$ by
\begin{align}
e_i' (P) = R \quad \mbox{and} \quad e_i'' (P) = Q.
\end{align}
Then we have the following commutation relations\,:
\begin{align}\label{eq:boson relation}
e_i' f_j = q_i^{-a_{ij}} f_j e_i' + \del_{ij} , \quad
e_i'' f_j = q_i^{a_{ij}} f_j e_i'' + \del_{ij}
\quad \mbox{ for any } i, j \in I.
\end{align}

Recall that $c_i \seteq -a_{ii}/2$ and set
\begin{align*}
e_i'{}^{(l)} = \begin{cases}
(e_i')^{l} \quad &\mbox{if }\ a_{ii}=2, \\
\dfrac{(e_i')^{l}}{\{ l \}_i!} &\mbox{if }\ a_{ii}\le 0.
\end{cases}
\end{align*}

By the commutation relation \eqref{eq:boson relation},
we have
\begin{align}\label{eq:come'f}
&e_i'{}^{(n)}f_j^{(m)}=
\begin{cases}
\sum\limits_{k=0}^n q_i^{-2nm+(n+m)k-k(k-1)/2}
{\genfrac{[}{]}{0pt}{0}{n}{k}}_i f_i^{(m-k)}
e_i'{}^{(n-k)}
&\mbox{if } i=j \mbox{ and }a_{ii}=2,\\
\sum\limits_{k=0}^m q_i^{-c_i(-2nm+(n+m)k-k(k-1)/2)}
{\genfrac{\{}{\}}{0pt}{0}{m}{k}}_i
f_i^{(m-k)}e_i'{}^{(n-k)}
&\mbox{if } i=j \mbox{ and } a_{ii}\le 0,\\
q_i^{-nma_{ij}}f_j^{(m)}e_i'{}^{(n)}
&\mbox{if } i \ne j.
\end{cases}
\end{align}

\smallskip

By a similar argument in \cite{Kas91},
one can show that there exists a unique non-degenerate
symmetric bilinear form on $\Um$ satisfying
\begin{align}
(1, 1) = 1, \quad (f_i P, Q) = (P, e_i' Q)
\quad \mbox{for all} \ \ i\in I.
\end{align}

Let $\star$ be the anti-involution of $\Uq$ defined by
\[ e_i^\star=e_i,\ f_i^\star=f_i \mbox{ and } (q^h)^\star=q^{-h}
\quad \mbox{for all} \ \ i\in I, h\in \pv. \]
Then we have
\begin{align}\label{eq:star}
(P^\star, Q^\star)=(P,Q) \quad \mbox{for any }P, Q \in \Um.
\end{align}

We define the following operator on $\Um$ :
\begin{align}\label{eq:proj}
P_i = \sum_{n \ge 0} (-1)^n q_i^{c_in (n-1)/2}
f_i^{(n)} e_i'{}^{(n)}.
\end{align}
Then $P_i$ is the projection operator to $\ker(e_i')$
with respect to the decomposition $\Um=\ker(e_i')\op f_i\Um$.
More precisely, we have the following result.
\begin{prop} \label{prop:decomp}
For each $i \in I$, every $u \in \Um$ can be uniquely expressed as
\begin{align} \label{eq:decomp}
u = \sum_{l \ge 0} f_i^{(l)} u_l ,
\end{align}
where $e_i' u_l = 0$ for every $l \ge 0$ and $u_l = 0$ for $l \gg 0$.
Moreover, we have
\[ u_{l} = q_i^{-c_i l (l-1)/2} P_i e_i'{}^{(l)} u. \]
\end{prop}

We call \eqref{eq:decomp} the {\em \str} of $u$.
\begin{defi}
For $i \in I$ and $u \in \Um$, let $u = \sum_{l \ge 0} f_i^{(l)} u_{l}$
be the \str.
Then the {\it Kashiwara operators} $\eit$ and $\fit$ on $\Um$ are defined
by
\begin{align}
\eit \, u = \sum_{l\ge 1} f_i^{(l-1)} u_l, \qquad
\fit \, u = \sum_{l\ge 0} f_i^{(l+1)} u_l.
\end{align}
\end{defi}

Note that the left multiplication operator
$f_i \cl \Um \to \Um$ is injective.
Hence we have
\begin{align}
\eit\circ\fit=\id_{\Um}.
\end{align}

\begin{defi}
A {\it crystal basis} of $\Um$ is a pair $(L, B)$, where
\begin{tenumerate}
\item $L$ is a free $\Az$-submodule of $\Um$ such that $\Um \simeq \Q(q)
\ot_{\Az} L$,
\item $B$ is a $\Q$-basis of $L/qL \simeq \Q \ot_{\Az} L$,
\item $\eit L \subset L$ and $\fit L \subset L$,
\item $\eit B \subset B \cup \ze$ and $\fit B \subset B$,
\item For $b, b' \in B$, we have $\fit b = b'$ if and only if $b = \eit
b'$.
\end{tenumerate}
\end{defi}

\begin{prop}\label{prop:l_0 for Um}
Let $(L, B)$ be a crystal basis of $\Um$.
For $i \in I$ and $u \in L_{-\al}$,
consider the \str
\begin{align*}
u = \sum_{l \ge 0} f_i^{(l)} u_l \quad \mbox{with} \quad e_i' u_l = 0 .
\end{align*}
Then the following statements ar true.
\begin{enumerate}
\item $u_l \in L$ for every $l \ge 0$.
\item If $u+qL \in B$, then there exists a non-negative integer
$l_0$ such that
\begin{enumerate}
\item $u_l \in q L$ for every $l \ne l_0$,
\item $u_{l_0} + qL \in B$,
\item $u \equiv f_i^{(l_0)} u_{l_0} \mod q L$.
\end{enumerate}
\end{enumerate}
\end{prop}

\medskip

\section{{\bf Existence of Crystal Bases}}

In this section,
we will prove the existence of crystal bases for
the algebra $\Um$ and for $\Uq$-modules in the category $\CO_{int}$.
Assuming the existence theorem,
the uniqueness of crystal bases can be proved by the same argument
as in \cite{Kas91} or \cite[Section 5.2]{HK2002}.

Let $\la \in P^+$ be a dominant integral weight
and let $\vla$ be the irreducible highest weight $\Uq$-module
with highest weight $\la$
and highest weight vector $v_\la$.
Let $\lla$ be the $\Az$-submodule of $\vla$
generated by the vectors of the form
\[ \wf 1 \dotsb \wf r v_\la \quad (r \ge 0, \ i_1,\dotsc, i_r \in I) \]
and set
\begin{align*}
\bla & \seteq \set{\wf 1 \dotsb \wf r v_\la + q \lla}
{r \ge 0, \ i_1, \dotsc, i_r \in I} \setminus \ze\\
&\subset \lla/q\lla.
\end{align*}
Similarly, we define $\linf$ to be
the free $\Az$-submodule of $\Um$ generated by the vectors of the form
\[ \wf 1 \dotsb \wf r \cdot 1 \quad (r \ge 0, \ i_1, \dotsc, i_r \in I) \]
and set
\[ \binf \seteq \set{\wf 1 \dotsb \wf r \cdot 1 + q \linf}
{r \ge 0, \ i_1, \dotsc, i_r \in I}\subset\linf/q\linf. \]

Our goal is to prove:
\begin{thm}\label{thm:existence}\hfill
\begin{enumerate}
\item The pair $(\lla, \bla)$ is a crystal basis of $\vla$.
\item The pair $(\linf, \binf)$ is a crystal basis of $\Um$.
\end{enumerate}
\end{thm}

It is easy to see that $\lla$ has the weight space decomposition
$\lla = \bop_{\mu \in P} \lla_\mu$,
where $\lla_\mu$ is the $\Az$-submodule of $\lla$ generated by
the vectors of the form $\wf 1 \dotsb \wf r v_\la$ such that
$\la - (\al_{i_1} + \dotsb + \al_{i_r}) = \mu$.
Similarly, $\bla$ has the weight set decomposition
$\bla = \bigsqcup_{\mu \in P} \bla_\mu$.
By the definition, we have, for $\al\in Q^+ \setminus \ze$
\[ \lla_{\la-\al}=\sum_{i \in I}\fit (\lla_{\la-\al+\al_i})\quad
\mbox{and}\quad
\bla_{\la-\al}=\Bigl(\bigcup_{i \in I}\fit (\bla_{\la-\al+\al_i})\Bigr)
\setminus \ze. \]
Since $\vla_\mu= \sum_i \fit \vla_{\mu+\al_i}$,
we can easily see that
$\vla \cong \Q(q) \ot_{\Az} \lla$.

Hence in order to prove (a),
it remains to prove the following statements\,:

$(1)$ $\eit \lla \subset \lla$ for every $i \in I$,

$(2)$ $\eit \bla \subset \bla \cup \ze $ for every $i \in I$,

$(3)$ for all $i \in I$ and $b, b' \in \bla$,
$\fit b = b'$ if and only if $b = \eit b'$,

$(4)$ $\bla$ is a $\Q$-basis of $\lla /q \lla$.

\smallskip

Similarly, to prove (b), we need to prove\,:

$(1)$ $\eit \linf \subset \linf$ for every $i \in I$,

$(2)$ $\eit \binf \subset \binf \cup \ze $ for every $i \in I$,

$(3)$ for all $i \in I$ and $b, b' \in \binf$,
$\fit b = b'$ if and only if $b = \eit b'$,

$(4)$ $\binf$ is a $\Q$-basis of $\linf /q \linf$.

We will prove these statements
using an interlocking induction argument on weights,
called the {\it grand-loop argument}, employed in \cite{Kas91}.

For $\la\in P^+$, let $\pi_\la\cl \Um\to \vla$
be the $\Um$-module homomorphism given by $P\mt P v_\la$.
For dominant integral weights $\la, \mu \in P^+$,
there exist unique $\Uq$-module homomorphisms
\begin{align*}
&\ophi \cl \vlm \to \vla \ot \vmu,\\
&\opsi \cl \vla \ot \vmu \to \vlm
\end{align*}
satisfying
\[ \ophi (v_{\la+\mu}) = v_\la \ot v_\mu,
\quad \opsi (v_\la \ot v_\mu) = v_{\la+\mu} . \]
It is clear that $\opsi \circ \ophi = \mbox{id}_{\vlm}$.

Let $(\ ,\ )$ denote the symmetric bilinear form
on irreducible highest weight $\Uq$-modules
defined in \eqref{eq:bilinear form},
and let us define a symmetric bilinear form on $\vla \ot \vmu$ by
\begin{align}
(u_1 \ot u_2 , v_1 \ot v_2) = (u_1, v_1) (u_2 ,v_2)
\end{align}
for $u_1, v_1 \in \vla$ and $u_2, v_2 \in \vmu$.
Then this form on $\vla \ot \vmu$ also satisfies
\begin{align}
(q^h u, v) = (u, q^h v), \quad
(f_i u, v) =(u, q_i^{-1} K_i e_i v), \quad
(e_i u, v) = (u, q_i f_i K_i^{-1} v)
\end{align}
for all $i \in I$ and $u, v \in \vla \ot \vmu$.
Moreover, we have
\begin{align}
(\opsi (u), v) = (u, \ophi (v))
\end{align}
for all $u \in \vla \ot \vmu$ and $v \in \vlm$.

Let $r \in \Z_{\ge 0}$ be a non-negative integer,
and set
\[Q_+(r) \seteq \set{\al \in Q_+ }{|\al| \le r } .\]
For $\la, \mu \in P^+$ and $\al \in Q_+(r)$,
we will prove that the following interlocking inductive statements are
true,
which would complete the proof of Theorem \ref{thm:existence}.

\begin{enumerate}\item[]\begin{enumerate}
\renewcommand{\labelenumii}{$\mathbf{\Alph{enumii}}(r)$ :}
\item $\eit \lla_{\la-\al} \subset \lla$ for every $i \in I$.

\item $\eit \bla_{\la-\al} \subset \bla \cup \ze $
for every $i \in I$.

\item For all $i \in I$, $b \in \bla_{\la-\al+\al_i}$
and $b' \in \bla_{\la-\al}$, we have
\[ \fit b = b' \mbox{ if and only if } b = \eit b' . \]

\item $\bla_{\la-\al}$ is a $\Q$-basis of
$\lla_{\la-\al}/q \lla_{\la-\al}$.

\item $\ophi (\llm_{\la+\mu-\al}) \subset \lom$.

\item $\opsi ((\lom)_{\la+\mu-\al}) \subset \llm$.

\item $\opsi ((\bla \ot \bmu)_{\la+\mu-\al}) \subset \blm \cup \ze$.

\item $\eit \linf_{-\al} \subset \linf$ for every $i \in I$.

\item $\eit \binf_{-\al} \subset \binf \cup \ze$
for every $i \in I$.

\item For all $i \in I$ and $b \in \binf_{-\al}$,
if $\eit b \ne 0$, then $b = \fit \eit b$.

\item $\binf_{-\al}$ is a $\Q$-basis of $\linf_{-\al}/q \linf_{-\al}$.

\item $\pi_\la (\linf_{-\al}) = \lla_{\la-\al}$.

\item For any $P \in \linf_{-\al+\al_i}$, we have
\[ \fit (P v_\la) \equiv (\fit P) v_\la \mod q \lla . \]

\item If $\ovl \pi_\la \cl (\linf / q \linf)_{-\al}
\longrightarrow (\lla / q \lla)_{\la-\al}$
is the induced homomorphism, we have
\[ \set{b \in \binf_{-\al}}{\ovl \pi_\la (b) \ne 0}
\isoto \bla_{\la-\al}. \]

\item If $b \in \binf_{-\al}$
satisfies $\ovl \pi_\la (b) \ne 0$, then we have
$\eit \ovl \pi_\la (b) = \ovl \pi_\la (\eit b)$.
\end{enumerate}\end{enumerate}

Note that the above statements are true for $r=0$ and $r=1$.
From now on, we assume that $r \ge 2$ and
the statements $\mathbf{A}(r-1), \dotsc, \mathbf{O}(r-1)$ are true.

\smallskip
In the sequel, for $\la\in P^+$,
``$\la \gg 0$'' means
``$\la(h_i) \gg 0$ for all $i \in I$''.
For instance,
``for any $P \in \linf$, we have $Pv_\la \in \lla$ for $\la \gg 0$''
is paraphrased as
``for any $P \in \linf$, there exists a positive integer $m$
such that $Pv_\la \in \lla$ whenever $\la \in P^+$ satisfies
$\la(h_i)\ge m$ for all $i \in I$''.

\begin{lem}\label{lem:conn}
Let $\al \in Q_+(r-1)$ and $b\in \bla_{\la-\al}$.
If $\eit b=0$ for every $i$, then we have
$\al=0$ and $b=v_\la$.
\end{lem}
\begin{proof}
If $\al \ne 0$, then there exist $i$ and $b'\in \bla_{\la-\al+\al_i}$
such that $b=\fit b'$.
Then $\mathbf{C}(r-1)$ implies
$\eit b=b' \ne 0$.
\end{proof}
\begin{lem}\label{lem:la-decomp}
For $\al \in Q_+(r-1)$,
let $u = \sum_{n=0}^N f_i^{(n)} u_n \in \vla_{\la-\al}$ be the \str,
where $u_n \in \ker e_i\cap \vla_{\la-\al+n \al_i}$.
\begin{enumerate}
\item If $u \in \lla$, then $u_n \in \lla$ for every $n \ge 0$.
\item If $u + q \lla \in \bla$, then there exists a non-negative integer
$n_0$ such that $u \equiv f_i^{(n_0)} u_{n_0} \mod q \lla$
and $u_n \in q \lla$ for $n \ne n_0$.
\end{enumerate}
\end{lem}
\begin{proof}
(a) We will use the induction on $N$.
If $N=0$, there is nothing to prove.
If $N > 0$, $\mathbf{A}(r-1)$ implies
$\eit u = \sum_{n=1}^N f_i^{(n-1)} u_n \in \lla$.
By the induction hypothesis, $u_n \in \lla$ for every $n = 1, \dotsc, N$.
Hence $u_0 = u - \sum_{n=1}^N \fit^{n} u_n \in \lla$,
which proves (a).

\smallskip

(b) If $N=0$, our assertion is trivial.
If $N > 0$, by $\mathbf{B}(r-1)$, we have
$\eit u + q \lla \in \bla \cup \ze$.

If $\eit u \in q \lla$,
then (a) implies that $u_n \in q \lla$ for every
$n \ge 1$, and $u \equiv u_0 \mod q \lla$.

If $\eit u \notin q \lla$, then $\eit u + q \lla \in \bla$.
By the induction hypothesis, there exists $n_0 \ge 1$ such that
$\eit u \equiv f_i^{(n_0 - 1)} u_{n_0} \mod q \lla$.
Hence by $\mathbf{C}(r-1)$, we obtain
\[ u \equiv \fit \eit u \equiv \fit f_i^{(n_0 - 1)} u_{n_0}
=f_i^{(n_0)} u_{n_0} \mod q \lla. \]
Therefore $u_n \in q\lla$ for $n \ne n_0$ by (a).
\end{proof}

\begin{lem}\label{lem:infty-decomp}
For $\al \in Q_+(r-1)$,
let $u = \sum_{n=0}^N f_i^{(n)} u_n \in \Um_{-\al}$
be the \str.
\begin{enumerate}
\item If $u \in \linf$, then $u_n \in \linf$ for every $n \ge 0$.
\item If $u + q \linf \in \binf$,
then there exists a non-negative integer $n_0$
such that $u \equiv f_i^{(n_0)} u_{n_0} \mod q \linf$
and $u_n \in q \linf$ for $n \ne n_0$.
\end{enumerate}
\end{lem}
\begin{proof}
The proof is similar to that of the preceding lemma.
\end{proof}

For $\al \in Q_+(r-1)$ and $b \in \bla_{\la-\al}$, we set
$\vep_i (b) = \max \set{n \ge 0}{\eit^n b \ne 0}$.
We define
\begin{align*}
\vp_i (b)& =
\begin{cases}
\vep_i(b)+\lan h_i,\wt(b)\ran &\mbox{if } i \in I^\re, \\
0 &\mbox{if } i \in I^\im \mbox{ and } \lan h_i, \wt(b)\ran = 0,\\
\infty &\mbox{if } i \in I^\im \mbox{ and } \lan h_i, \wt(b) \ran > 0.
\end{cases}
\end{align*}

Note that $\vp_i(b)=0$ implies $\fit b=0$.

\begin{lem}\label{lem:long}
Suppose $\al, \be \in Q_+(r-1)$ and
$b \in \bla_{\la-\al}$, $b' \in \bmu_{\mu-\be}$.
\begin{enumerate}
\item For each $i \in I$, we have
\begin{align*}
& \fit (\lla_{\la-\al} \ot \lmu_{\mu-\be}) \subset \lom, \\
& \eit (\lla_{\la-\al} \ot \lmu_{\mu-\be}) \subset \lom.
\end{align*}
\item For each $i \in I$, the tensor product rule holds in $(\lom) / q
(\lom)$\,{\rm :}
\begin{align*}
&\fit (b \ot b') =
\begin{cases} \ \fit b \ot b' &\mbox{if } \vp_i (b) > \vep_i (b'),\\
\ b \ot \fit b' &\mbox{if } \vp_i (b) \le \vep_i (b'),
\end{cases} \\
&\eit (b \ot b') =
\begin{cases} \ \eit b \ot b' &\mbox{if } \vp_i (b) \ge \vep_i (b'),\\
\ b \ot \eit b' &\mbox{if } \vp_i (b) < \vep_i (b').
\end{cases}
\end{align*}
\item For each $i \in I$, we have
\begin{align*}
& \fit (\bla_{\la-\al} \ot \bmu_{\mu-\be})
\subset (\bla \ot \bmu) \cup \ze, \\
& \eit (\bla_{\la-\al} \ot \bmu_{\mu-\be})
\subset (\bla \ot \bmu) \cup \ze.
\end{align*}
\item If $\eit (b \ot b') \ne 0$,
then $b \ot b' \equiv \fit \eit (b \ot b')$.
\item If $\eit (b \ot b') =0$ for every $i \in I$, then $b = v_\la$.
\item For each $i \in I$, we have
\[ \fit (b \ot v_\mu) \equiv \fit b \ot v_\mu
\quad \mbox{\rm{or}} \quad \fit b = 0. \]
\item For any sequence of indices $i_1, \dotsc, i_r \in I$, we have
\[ \wf 1 \dotsb \wf r (v_\la \ot v_\mu)
\equiv (\wf 1 \dotsb \wf r v_\la) \ot v_\mu \mod q \lom \]
or \hskip 0.3 cm $\wf 1 \dotsb \wf r v_\la \in q \lla$.
\end{enumerate}
\end{lem}
\begin{proof}
(a) By Lemma \ref{lem:la-decomp}(a),
it suffices to prove the following statement\,:
for $u \in \ker(e_i)\cap\lla_{\la-\al+n \al_i}$ and
$v \in\ker(e_i)\cap \lmu_{\mu-\be+m \al_i}$, we have
\begin{align*}
&\fit (f_i^{(n)} u \ot f_i^{(m)} v ) \in \lom.,\\
&\eit (f_i^{(n)} u \ot f_i^{(m)} v ) \in \lom.
\end{align*}
Let $L$ be the free $\Az$-submodule of $\vla \ot \vmu$
generated by the vectors
$f_i^{(s)} u \ot f_i^{(t)} v$ $(s, t \ge 0)$.
Then by the tensor product rule for $U_i$-modules,
we have $\eit L \subset L$ and $\fit L \subset L$.
Hence our assertion follows immediately
from the fact that $L \subset \lom$.

\noindent
(b), (d) By Lemma \ref{lem:la-decomp}(b),
we may assume that $b = f_i^{(n)} u + q \lla$
and $b' = f_i^{(m)} v + q \lmu$ with $e_iu=0$ and $e_iv=0$.
Let $L$ be the $\Az$-submodule generated by
$f_i^{(s)} u \ot f_i^{(t)} v$ $(s, t \ge 0)$.
Then by the tensor product rule for $U_i$-modules,
our assertion holds in $L / q L$.
Since $L \subset \lom$, we are done.

\noindent
The statements (c) is a simple consequence of (b).

\noindent
(e) By (b), we have $\eit b=0$
and Lemma \ref{lem:conn} implies the desired result.

\noindent
(f) By the tensor product rule (b),
we have $\fit (b \ot v_\mu) \equiv \fit b \ot v_\mu$
unless $\vp_i (b) = 0$, in which case $\fit b = 0$.

\noindent
(g) This is an immediate consequence of (f).
\end{proof}

\begin{prop}[$\mathbf{E}(r)$]\label{prop:e(r)}
For every $\al \in Q_+(r)$, we have
\[ \ophi (\llm_{\la+\mu-\al}) \subset \lom. \]
\end{prop}
\begin{proof}
Recall that
$\llm_{\la+\mu-\al} = \sum_i \fit \llm_{\la+\mu-\al+\al_i}$.
By $\mathbf{E}(r-1)$ and Lemma \ref{lem:long}(a),
we obtain the desired result.
\end{proof}

\begin{lem}\label{lem:key-lemma}
Let $i_1, \dotsc, i_r$ be a sequence of indices in $I$.
Suppose that $i_t \ne i_{t+1} = \dotsb = i_r$
for a positive integer $t<r$.
Then for any dominant integral weight
$\mu \in P^+$ and $\la = \La_{i_t}$, we have
\[ \wf 1 \dotsb \wf r (v_\la \ot v_\mu) \equiv b \ot b'
\quad \mod q \lom \]
for some $b \in \bla_{\la-\al}\cup\ze$,
$b' \in \bmu_{\mu-\be} \cup \ze$ and
$\al, \be \in Q_+(r-1) \setminus \ze$.
\end{lem}
\begin{proof}
Since $\La_{i_t} (h_{i_r}) = 0$, we have $f_{i_r} v_\la = 0$.
Hence we get
\begin{align*}
f_{i_r} (v_\la \ot v)
= (f_{i_r} \ot 1 + K_{i_r} \ot f_{i_r}) (v_\la \ot v)
= v_\la \ot f_i v
\end{align*}
for every $v \in \vmu$,
which gives
\begin{align*}
\wf {t+1} \dotsb \wf r (v_\la \ot v_\mu)
= f_{i_r}^{(r-t)} (v_\la \ot v_\mu)
= v_\la \ot f_{i_r}^{(r-t)} v_\mu .
\end{align*}
Since $e_{i_t} (v_\la \ot f_{i_r}^{(r-t)} v_\mu) = 0$, we have
\begin{align*}
\wf t \dotsb \wf r (v_\la \ot v_\mu)
&= \wf t (v_\la \ot f_{i_r}^{(r-t)} v_\mu)
= f_{i_t} (v_\la \ot f_{i_r}^{(r-t)} v_\mu) \\
&= (f_{i_t} \ot 1 + K_{i_t} \ot f_{i_t})
(v_\la \ot f_{i_r}^{(r-t)} v_\mu) \\
&= f_{i_t} v_\la \ot f_{i_r}^{(r-t)} v_\mu
+ q_{i_t}^{\la(h_{i_t})} v_\la \ot f_{i_t} f_{i_r}^{(r-t)} v_\mu \\
&= \wf t v_\la \ot \wf r^{r-t} v_\mu
+ q_{i_t} v_\la \ot \wf t \wf r^{r-t} v_\mu \\
&\equiv \wf t v_\la \ot \wf r^{r-t} v_\mu \mod q \lom.
\end{align*}
Hence Lemma \ref{lem:long}(b) yields our assertion.
\end{proof}

\begin{lem}\label{lem:phi-1}
For every $\al \in Q_+(r)$, we have
\[ (\lom)_{\la+\mu-\al} =
\sum_i \fit (\lom)_{\la+\mu-\al+\al_i}
+ v_\la \ot \lmu_{\mu-\al}. \]
\end{lem}
\begin{proof}
Let $L$ be the left-hand side and
$L'$ be the right-hand side of the above equation.
We already know $L' \subset L$ by Lemma \ref{lem:long}(a).
By Lemma \ref{lem:long}(e),
for any $\be \in Q_+(r-1) \setminus \ze$ and
$b \ot b' \in \bla_{\la-\be} \ot \bmu_{\mu-\al+\be}$,
there exists $i \in I$ such that $\eit (b \ot b') \ne 0$.
By Lemma \ref{lem:long}(d), we have
\[ b \ot b' \equiv \fit \eit (b \ot b') \mod q \lom. \]
Hence $\lla_{\la-\be} \ot \lmu_{\mu-\al+\be} \subset L' + q L$,
and we obtain
\begin{align*}
L &= \lla_{\la-\al} \ot v_\mu
+ \sum_{\be \in Q_+(r-1) \setminus \ze}
\lla_{\la-\be} \ot \lmu_{\mu-\al+\be}
+ v_\la \ot \lmu_{\mu-\al} \\
& \subset \lla_{\la-\al} \ot v_\mu + L' + qL .
\end{align*}
Note that for any $\wf 1 \dotsb \wf r v_\la \in \bla_{\la-\al}$,
Lemma \ref{lem:long}(f) yields
\[ (\wf 1 \dotsb \wf r v_\la) \ot v_\mu \equiv
\wf 1 \wf 2 \dotsb \wf r (v_\la \ot v_\mu)
\mod q \lom. \]
It follows that $\lla_{\la-\al} \ot v_\mu \subset L' + q L$
and hence $L\subset L' + q L$.
Then by Nakayama's Lemma, we conclude that $L = L'$.
\end{proof}

For $\la, \mu \in P^+$, define a linear map
$\os \cl \vla \ot \vmu \to \vla$ by
\begin{align}\label{eq:s}
\begin{aligned}
& \os (u \ot v_\mu) = u \quad \mbox{for every } u \in \vla,\\
& \os \bigl( \vla \ot \sum_i f_i \vmu \bigr) = 0.
\end{aligned}
\end{align}
Note that $\os$ is a $\Um$-module homomorphism.

\begin{lem}\label{lem:s}\hfill
\begin{enumerate}
\item For all $\la, \mu \in P^+$, we have $\os (\lom) = \lla$.
\item For all $\al \in Q_+(r-1)$
and $w \in (\lom)_{\la+\mu-\al}$, we have
\[ \os \circ \fit (w)
\equiv \fit \circ \os (w) \quad \mod q \lla. \]
\end{enumerate}
\end{lem}
\begin{proof}
(a) is obvious.

\noindent
(b) By Lemma \ref{lem:la-decomp},
we may assume that $w = f_i^{(n)} u \ot f_i^{(m)} v$
with $u \in \lla$, $v \in \lmu$ and $e_i u = e_i v =0$.
Let $L$ be the $\Az$-submodule generated by
$f_i^{(s)} u \ot f_i^{(t)} v$ $(s, t \ge 0)$.
Then the tensor product rule for $U_i$-modules yields
\[ \fit w \equiv f_i^{(n+1)} u \ot f_i^{(m)} v
\quad \mbox{or} \quad f_i^{(n)} u \ot f_i^{(m+1)} v \mod q L. \]
Observe that
$\os (\fit w) \equiv\fit (\os w) \equiv 0 \mod q L$
unless $m=0$ and $v \in \lmu_\mu = \Az v_\mu$.
If $m=0$ and $v = v_\mu$, then we have
\begin{align*}
\fit w = \fit (f_i^{(n)} u \ot v_\mu) \equiv
\begin{cases}
f_i^{(n+1)} u \ot v_\mu \mod q L
\quad & \mbox{if } f_i^{(n+1)} u \ne 0, \\
f_i^{(n)} u \ot f_i v_\mu \mod q L
& \mbox{if } f_i^{(n+1)} u = 0 .
\end{cases}
\end{align*}
In both cases, we have
$\os (\fit w) \equiv f_i^{(n+1)} u = \fit f_i^{(n)} u
= \fit \os (w) \mod q \lla$.
\end{proof}

\begin{lem}\label{lem:pre-m(r)}
Let $\al\in Q_+(r)$ and $m\in \Z_{\ge 0}$.
Then there exists a positive integer $N$ satisfying
the following properties$:$
if $\la\in P^+$ and $u\in \vla_{\la-\al}$ satisfies
$q_i^{-1}K_i e_iu\in q^N \lla$,
then we have
\begin{enumerate}
\item
$\fit^mu\equiv f_i^{(m)}u\mod q \lla$,
\item
$\eit f_i^{(m)}u \equiv f_i^{(m-1)}u \mod q\lla$
if $i \in I^\re$ and $m\le \lan h_i,\la-\al \ran$
or if $i \in I^\im$ and $\la(h_i)>0$.
\end{enumerate}
\end{lem}
\begin{proof}
Let $u =\sum_{n=0}^r f_i^{(n)}u_n$ be the \str\ of $u$.
Then we have
\[ q_i^{-1}K_i e_i u =\sum_{n=1}^r q_i^{-1} K_i e_i f_i^{(n)}u_n \]
with $q_i^{-1} K_i e_i f_i^{(n)}u_n \in f_i^{(n-1)} (\ker e_i)$.
Hence Lemma \ref{lem:la-decomp} implies that
$q_i^{-1}K_i e_i f_i^{(n)}u_n\in q^{N}\lla$ for $n\ge 1$.
By \eqref{eq:qKe}, we have
$f_i^{(n-1)}u_n \in \Az \bigr(q_i^{-1} K_i e_i f_i^{(n)}u_n \bigr)$
for $n \ge 1$.
We then conclude
$u_n \in q^{N}\lla$ for $n \ge 1$
by Lemma \ref{lem:la-decomp}.

Hence if $N \gg 0$, the elements
$f_i^{(m)}f_i^{(n)}u_n$,
$f_i^{(m-1)}f_i^{(n)}u_n$, $f_i^{(m+n)}u_n$
belong to $q \lla$ for $n \ge 1$, which implies
\[ \fit^m u = \sum_{n=0}^r f_i^{(n+m)}u_n
\equiv f_i^{(m)}u_0 \mod q \lla \]
and
\[ f_i^{(m)}u = \sum_{n=0}^r f_i^{(m)}f_i^{(n)}u_n
\equiv f_i^{(m)}u_0 \mod q \lla. \]
Hence we obtain
$\fit^m u\equiv f_i^{(m)}u$.
Similarly, we have
$\fit^{m-1} u = f_i^{(m-1)}u \equiv f_i^{(m-1)}u_0$.

On the other hand, we have
$f_i^{(m)}u =\sum_{n=0}^r a_nf_i^{(m+n)}u_n$,
where $a_n$ is a $q$-binomial coefficient or $1$
according as $i$ is real or imaginary.
Since
$a_n f_i^{(n+m-1)}u_n \in a_n q^N \lla \subset q\lla$
for $n \ge 1$ and $N \gg 0$, \eqref{eq:te} yields
\begin{align*}
\eit f_i^{(m)}u
&=\sum_{f_i^{(m+n)}u_n \ne 0} a_nf_i^{(n+m-1)}u_n\\
&\equiv
\begin{cases}
f_i^{(m-1)}u_0 &\mbox{if } f_i^{(m)}u_0 \ne 0 \mod q \lla, \\
0 &\mbox{if } f_i^{(m)}u_0=0 \mod q \lla .
\end{cases}
\end{align*}
However, by the conditions in (b),
$f_i^{(m)}u_0=0$ implies $u_0=0$.
Therefore we conclude $\eit f_i^{(m)}u \equiv f_i^{(m-1)}u_0
\equiv f_i^{(m-1)}u \mod q\lla$.
\end{proof}

\begin{lem}\label{lem:infty-la}
Let $\al \in Q_+(r)$ and $P \in \Um_{-\al}$.
Then for every $\la \gg 0$, we have
\begin{align*}
\begin{cases}
(\fit P) v_\la \equiv \fit (P v_\la) & \mod q \lla, \\
(\eit P) v_\la \equiv \eit (P v_\la) & \mod q \lla.
\end{cases}
\end{align*}
\end{lem}
\begin{proof}
We may assume that
$P=f_i^{(n)}Q$ with $e_i'Q=0$.
Then
\begin{align*}
q_i^{-1}K_i e_i(Qu_\la)
&=(1-q_i^2)^{-1}(e'_iQ-K_i^2e_i'' Q)u_\la\\
&=-(1-q_i^2)^{-1}q_i^{2 \lan h_i,\la+\al_i+\wt(Q)\ran}(e_i'' Q)u_\la
\end{align*}
Hence by the preceding lemma, we have for $\la \gg 0$
\[ \fit(f_i^{(n)}Q v_\la) \equiv f_i^{(n+1)}Q v_\la
\quad \mbox{ and } \quad
\eit(f_i^{(n)}Q v_\la)\equiv f_i^{(n-1)}Q v_\la. \]
\end{proof}

\begin{prop}[$\mathbf{M}(r)$]\label{prop:m(r)}
For any $\la \in P^+$, $\al \in Q_+(r-1)$ and $P \in \linf_{-\al}$,
we have
\begin{align*}
(\fit P) v_\la \equiv \fit (P v_\la) & \mod q \lla.
\end{align*}
In particular, we have
$(\wf 1 \dotsb \wf r \cdot 1) v_\la
\equiv \wf 1 \dotsb \wf r v_\la \mod q \lla$.
\end{prop}
\begin{proof}
Take $\mu \gg 0$.
Since $\la + \mu \gg 0$, Lemma \ref{lem:infty-la} yields
\begin{align*}
(\fit P) v_{\la+\mu} \equiv \fit (P v_{\la+\mu}) \mod q \llm .
\end{align*}
By Proposition \ref{prop:e(r)}, applying $\ophi$ gives
\begin{align}\label{eq:fit P}
(\fit P) (v_\la \ot v_\mu) \equiv \fit (P (v_\la \ot v_\mu)) \mod q \lom .
\end{align}

On the other hand, we have
\[ P(v_\la \ot v_\mu) = \ophi (Pv_{\la+\mu}) \in \lom \]
by $\mathbf{L}(r-1)$ and $\mathbf{E}(r-1)$.
Applying $\os$ to \eqref{eq:fit P}, Lemma \ref{lem:s} implies
\begin{align*}
(\fit P) v_\la \equiv \fit (P v_\la) \mod q \lla ,
\end{align*}
as desired.
\end{proof}

\begin{prop}[$\mathbf{L}(r)$]\label{prop:l(r)}
For any $\la \in P^+$ and $\al \in Q_+(r)$, we have
\[ \pi_\la (\linf_{-\al}) = \lla_{\la-\al} . \]
\end{prop}
\begin{proof}
By Proposition \ref{prop:m(r)}, we have
\[ \pi_\la (\wf 1 \dotsb \wf r \cdot 1)
\equiv \wf 1 \dotsb \wf r v_\la \mod q \lla. \]
It follows that
\[ \pi_\la (\linf_{-\al})\subset \lla_{\la-\al} \mbox{ and }
\lla_{\la-\al} \subset \pi_\la (\linf_{-\al}) + q \lla_{\la-\al}. \]
Now the desired result follows from Nakayama's Lemma.
\end{proof}

\begin{coro}\label{coro:large-isom}
Let $\ovl \pi_\la \cl (\linf / q \linf)_{-\al}
\longrightarrow (\lla / q \lla)_{\la-\al}$
be the surjective $\Q$-linear map induced by $\pi_\la$.
\begin{enumerate}
\item
For $\be\in Q_+(r-1)$ and $b\in \binf_{-\be}$, we have
\[ \ovl \pi_\la(\fit b)=\fit\ovl \pi_\la(b). \]
\item For $\al \in Q_+(r)$ and $\la \in P^+$, we have
\begin{align*}
(\ovl \pi_\la \binf_{-\al}) \setminus \ze = \bla_{\la-\al} .
\end{align*}
\item Given $\al \in Q_+(r)$, take $\la \gg 0$.
Then $\pi_\la$ induces the isomorphisms
\[ \linf_{-\al} \isoto \lla_{\la-\al}, \quad
\binf_{-\al} \setminus \ze \isoto \bla_{\la-\al}. \]
\end{enumerate}
\end{coro}
\begin{proof}
(a) is an immediate consequence of Proposition \ref{prop:m(r)},
and (b) follows from (a).
Finally, (c) follows from the fact that
$\pi_\la\cl\Um_{-\al} \to \vla_{\la-\al}$ is an isomorphism
for $\la \gg 0$.
\end{proof}

Now we will prove the statements $\mathbf{A}(r)$ and $\mathbf{H}(r)$.
Fix $\la \in P^+$, $i \in I$
and $\al = \al_{i_1} + \dotsb + \al_{i_r} \in Q_+(r)$.
Let $T$ be a finite subset of $P^+$ containing $\la$ and
the fundamental weights $\La_{i_1}, \dotsc, \La_{i_r}$.
Since $T$ is a finite set, we can take a sufficiently large
$N_1 \ge 0$ such that
\[ \eit \lmu_{\mu-\al} \subset q^{-N_1} \lmu
\quad \mbox{for every} \quad \mu \in T. \]
Let $N_2 \ge 0$ be a non-negative integer such that
\[ \eit \linf_{-\al} \subset q^{-N_2} \linf. \]
Then for every $\mu \gg 0$, by Lemma \ref{lem:infty-la} and
Proposition \ref{prop:l(r)}, we have
\begin{align*}
\eit \lmu_{\mu-\al}
&= \eit ( \linf_{-\al} v_\mu )
\subset (\eit \linf_{-\al}) v_\mu + q \lmu_{\mu-\al} \\
& \subset q^{-N_2} \linf_{-\al} v_\mu + q \lmu_{\mu-\al} \subset q^{-N_2} \lmu.
\end{align*}
Therefore, given $\al\in Q_+(r)$,
there exists a non-negative integer $N \ge 0$ such that
\begin{align}\label{eq:small loop}
\begin{aligned}
& \eit \lmu_{\mu-\al} \subset q^{-N} \lmu &&\mbox{for every}
\quad \mu \gg 0,\\
& \eit \ltau_{\tau-\al} \subset q^{-N} \ltau &&\mbox{for every}
\quad \tau \in T, \\
& \eit \linf_{-\al} \subset q^{-N} \linf.
\end{aligned}
\end{align}

\begin{lem}\label{lem:uniform-n}
Given $\al \in Q_+(r)$,
let $N \ge 0$ be a non-negative integer satisfying \eqref{eq:small loop}.
Then for every $\mu \gg 0$ and $\tau \in T$, we have
\[ \eit \bigl((\ltau \ot \lmu)_{\tau+\mu-\al}\bigr) \subset q^{-N}
(\ltau \ot \lmu). \]
\end{lem}
\begin{proof}
For $u \in \ltau_{\tau-\be}$ and $v \in \lmu_{\mu-\ga}$
with $\al = \be + \ga$, let us show
$\eit(u \ot v)\in q^{-N}(\ltau \ot \lmu)$.
If $\be \ne 0$ and $\ga \ne 0$,
we already proved our assertion in Lemma \ref{lem:long}(a).

If $\be=0$, then $\ga=\al$ and we may assume $u = v_\tau$.
Let $v = \sum_m f_i^{(m)} v_m$ with $e_i v_m =0$ be the \str\ of $v$.
By \eqref{eq:small loop}, we have
\[ \eit v = \sum_{m\ge 1} f_i^{(m-1)} v_m \in q^{-N} \lmu. \]
Lemma \ref{lem:la-decomp} implies
 that $v_m \in q^{-N} \lmu$ for every $m \ge 1$.
Let $L$ be the $\Az$-submodule generated by the vectors
$f_i^{(s)} v_\tau \ot f_i^{(t)} v_m$ $(s,\ t \ge 0,\ m\ge 1)$.
Then we have $\eit L \subset L$.
It follows that
\[ \eit (v_\tau \ot v)
= \sum_{m \ge 1} \eit (v_\tau \ot f_i^{(m)} v_m) \in L
\subset q^{-N} \ltau \ot \lmu. \]

If $\be=\al$ and $\ga=0$,
a similar argument yields
$\eit (u \ot v_\mu) \in q^{-N} \ltau \ot \lmu$.
\end{proof}

\begin{lem}\label{lem:-n+1}
Given $\al \in Q_+(r)$, let $N \ge 1$ be a positive integer satisfying
\eqref{eq:small loop}.
Then we have
\begin{enumerate}
\item $\eit \lmu_{\mu-\al} \subset q^{-N+1} \lmu$ for all $\mu \gg 0$,
\item $\eit \ltau_{\tau-\al} \subset q^{-N+1} \ltau$
for every $\tau \in T$,
\item $\eit \linf_{-\al} \subset q^{-N+1} \linf$.
\end{enumerate}
\end{lem}
\begin{proof}
(a) Let $u = \wf 1 \dotsb \wf r v_\mu \in \lmu_{\mu-\al}$.
Suppose $i_1 = \dotsb = i_r$.
If $i=i_1$, then $\eit u = \eit f_i^{(r)} v_\mu = f_i^{(r-1)} v_\mu \in
\lmu$,
and if $i \ne i_1$, $\eit u = 0$.
Hence we may assume that there exists a positive integer $t<r$ such that
$i_t \ne i_{t+1} = \dotsb = i_r$.

Suppose $\mu \gg 0$ and set $\la_0 = \La_{i_t}$,
$\mu' = \mu - \la_0 \gg 0$.
By Lemma \ref{lem:key-lemma}, we have
\[ w \seteq \wf 1 \dotsb \wf r (v_{\la_0} \ot v_{\mu'})
\equiv v \ot v' \mod q L(\la_0) \ot L(\mu') \]
for some $v \in L(\la_0)_{\la_0-\be}$, $v' \in L(\mu')_{\mu'-\ga}$,
$\be, \ga \in Q_+(r-1) \setminus \ze$, $\al = \be + \ga$.
Hence Lemma \ref{lem:long}(a) and Lemma \ref{lem:uniform-n} yield
\begin{align*}
\eit w &\in L(\la_0) \ot L(\mu')
+ q \eit (L(\la_0) \ot L(\mu'))_{\la_0 +\mu' -\al} \\
&\subset L(\la_0) \ot L(\mu') + q^{-N+1} L(\la_0) \ot L(\mu')
= q^{-N+1} L(\la_0) \ot L(\mu') .
\end{align*}
Namely, we have
\[ \eit w \in q^{-N+1} (L(\la_0) \ot L(\mu'))_{\la_0+\mu'-\al+\al_i} . \]
Therefore, applying $\Psi_{\la_0, \mu'}$,
the statement $\mathbf{F}(r-1)$ yields
\[ \eit u = \eit \wf 1 \dotsb \wf r v_\mu \in q^{-N+1} \lmu . \]

(b) Let $\tau \in T$ and $u = \wf 1 \dotsb \wf r v_\tau \in
\ltau_{\tau-\al}$.
If $u \in q \ltau$, then our assertion follows immediately from
\eqref{eq:small loop}.
If $u \notin q \ltau$, then for any $\mu \in P^+$,
Lemma \ref{lem:long}(g) gives
\begin{align}\label{eq:former}
\wf 1 \dotsb \wf r (v_\tau \ot v_\mu)
\equiv u \ot v_\mu \mod q \ltau \ot \lmu .
\end{align}

If $\mu \gg 0$, then (a) implies
\[ \eit \wf 1 \dotsb \wf r v_{\tau+\mu} \in q^{-N+1} L(\tau+\mu) .\]
Applying $\Phi_{\tau, \mu}$, the statement $\mathbf{E}(r-1)$ gives
\[ \eit \wf 1 \dotsb \wf r (v_\tau \ot v_\mu)
\in q^{-N+1} \ltau \ot \lmu . \]
Along with \eqref{eq:former}, Lemma \ref{lem:uniform-n} yields
\begin{align} \label{eq:-n+1}
\begin{aligned}
\eit (u \ot v_\mu)
&\in q^{-N+1} \ltau \ot \lmu + q \eit (\ltau \ot \lmu) \\
&\subset q^{-N+1} \ltau \ot \lmu .
\end{aligned}
\end{align}
Let $u = \sum\limits_{m\ge 0} f_i^{(m)} u_m$ be the \str\ of $u$.
Then $\eit u = \sum\limits_{m\ge 1} f_i^{(m-1)} u_m \in q^{-N} \ltau$
implies that $u_m \in q^{-N} \ltau$ for $m \ge 1$.
Let $L$ be the $\Az$-submodule generated by
$f_i^{(s)} v_m \ot f_i^{(t)} v_\mu$ $(s, t \ge 0, m \ge 1)$.
Then $L \subset q^{-N} \ltau \ot \lmu$
and by the tensor product rule for $U_i$-modules,
we obtain
\begin{align*}
\eit (u \ot v_\mu) & = \sum_{m \ge 1} \eit (f_i^{(m)} u_m \ot v_\mu) \\
&\equiv \sum_{m \ge 1} f_i^{(m-1)} u_m \ot v_\mu
= \eit u \ot v_\mu \mod q L .
\end{align*}
Hence we have $\eit (u \ot v_\mu)\equiv\eit u \ot v_\mu
\mod q^{1-N}\ltau\ot \lmu$.
Together with \eqref{eq:-n+1}, we obtain
$\eit u \ot v_\mu\in q^{-N+1} \ltau \ot \lmu$, which implies
the desired result $\eit u \in q^{-N+1} \ltau$.

\smallskip

\noindent
(c) Let $P \in \linf_{-\al}$ and take $\mu \gg 0$.
By Lemma \ref{lem:infty-la}, we have
\begin{align*}
(\eit P) v_\mu \equiv \eit (P v_\mu) \mod q \lmu.
\end{align*}
By Proposition \ref{prop:l(r)}, we get
\[ \eit (P v_\mu) \in \eit \lmu_{\mu-\al} \subset q^{-N+1} \lmu, \]
which implies $(\eit P) v_\mu \in q^{-N+1} \lmu$.
Hence by Corollary \ref{coro:large-isom} (c),
we obtain $\eit P \in q^{-N+1} \linf$.
\end{proof}

\begin{prop}[$\mathbf{A}(r)$, $\mathbf{H}(r)$]\label{prop:a(r),h(r)}
For all $\la \in P^+$ and $\al \in Q_+(r)$, we have
\begin{align*}
\eit \linf_{-\al} \subset \linf, \quad \quad
\eit \lla_{\la-\al} \subset \lla.
\end{align*}
\end{prop}
\begin{proof}
By applying Lemma \ref{lem:-n+1} repeatedly, we conclude that $N = 0$,
and our assertion follows immediately.
\end{proof}

\begin{coro}\label{coro:non-zero}
For $\al \in Q_+(r)$, $\binf_{-\al}$ does not contain $0$.
\end{coro}
\begin{proof}
For $b \in \binf_{-\al}$, there exists $i\in I$ and
$b' \in \binf_{-\al + \al_i}$ such that
$b =\fit b'$.
Then $b' = \eit b$ does not vanish by $\mathbf{K}(r-1)$,
and hence $b \ne 0$.
\end{proof}

\begin{lem}\label{lem:e-inf-la}
For $\al\in Q_+(r)$, $i \in I$,
$\la \gg 0$ and $b\in\binf_{-\al}$, we have
\[ \ovl\pi_\la(\eit b)=\eit\ovl\pi_\la(b). \]
\end{lem}
\begin{proof}
This follows immediately from Lemma \ref{lem:infty-la}.
\end{proof}

Since Lemma \ref{lem:long} depends only on $\mathbf{A}(r-1)$, we have:

\begin{coro}\label{coro:long}
Let $\la, \mu \in P^+$ and $\al, \be \in Q_+(r)$.
\begin{enumerate}
\item If $u = \sum_{n\ge 0}f_i^{(n)} u_n \in \lla_{\la-\al}$
is the \str\ of $u$,
then $u_n \in \lla$ for every $n \ge 0$.
\item For each $i \in I$, we have
\begin{align*}
&\fit (\lla_{\la-\al} \ot \lmu_{\mu-\be}) \subset \lom, \\
&\eit (\lla_{\la-\al} \ot \lmu_{\mu-\be}) \subset \lom.
\end{align*}
\end{enumerate}
\end{coro}

In order to prove the statements $\mathbf{B}(r)$ and $\mathbf{O}(r)$,
we first prove:
\begin{lem}\label{lem:eit}
Let $\la$, $\mu \in P^+$ and $\al \in Q_+(r)$.
Then for every $u \in \lla_{\la-\al}$, we have
\[ \eit (u \ot v_\mu) \equiv (\eit u) \ot v_\mu \mod q \lom. \]
\end{lem}
\begin{proof}
By Corollary \ref{coro:long}(a),
we may assume $u = f_i^{(n)} w \ne 0$ with $e_i w =0$
and $w \in \lla_{\la-\al + n\al_i}$.
Let $L$ be the $\Az$-submodule generated by
$f_i^{(s)} w \ot f_i^{(t)} v_\mu$ $(s, t \ge 0)$.
Then the tensor product rule for $U_i$-modules gives
\[ \eit (f_i^{(n)} w \ot v_\mu )
\equiv f_i^{(n-1)} w \ot v_\mu \mod q L. \]
Our assertion immediately follows from
$L \subset \lom$.
\end{proof}

\begin{prop}[$\mathbf{I}(r)$]\label{prop:i(r)}
For every $\al \in Q_+(r)$, we have
\[ \eit \binf_{-\al} \subset \binf \cup \ze. \]
\end{prop}
\begin{proof}
Let $P = \wf 1 \dotsb \wf r \cdot 1 \in \linf_{-\al}$ and
$b = \wf 1 \dotsb \wf r \cdot 1 + q \linf \in \binf_{-\al}$.
If $i_1 = \dotsb = i_r$,
our assertion is true as we have seen
in the proof of Lemma \ref{lem:-n+1} (a).

Hence we may assume that there exists a positive integer $t<r$ such that
$i_t \ne i_{t+1} = \dotsb = i_r$.
Take $\mu \gg 0$ and set $\la_0 = \La_{i_t}$,  $\la=\la_0+\mu$.
Then Lemma \ref{lem:key-lemma} yields
\[ \wf 1 \dotsb \wf r (v_{\la_0} \ot v_\mu )
\equiv v \ot v' \mod q L(\la_0) \ot \lmu \]
for some $v \in L(\la_0)_{\la_0-\be}$,
$v' \in \lmu_{\mu-\ga}$ and
$\be$, $\ga \in Q_+(r-1) \setminus \ze$ with $\al = \be + \ga$
such that $v + q L(\la_0) \in B(\la_0)\cup\ze$
and $v' + q \lmu \in \bmu \cup \ze$.
Then the tensor product rule for $U_i$-modules gives
\begin{align*}
\eit \wf 1 \dotsb \wf r (v_{\la_0} \ot v_\mu)
&\equiv \eit (v \ot v') \\
&\equiv \eit v \ot v' \mbox{ or } v \ot \eit v' \mod q L(\la_0) \ot \lmu.
\end{align*}
By $\mathbf{B}(r-1)$, we have
\[ \eit \wf 1 \dotsb \wf r (v_{\la_0} \ot v_\mu) + q L(\la_0) \ot \lmu
\in (B(\la_0) \ot \bmu) \cup \ze. \]
Hence, by applying $\Psi_{\la_0,\mu}$, the statement $\mathbf{G}(r-1)$
gives
\[ \eit\ovl\pi_\la(b)
= \eit \wf 1 \dotsb \wf r v_\la + q \lla \in \bla \cup \ze . \]
Therefore, by Lemma \ref{lem:e-inf-la} and Corollary
\ref{coro:large-isom}(c),
we conclude $\eit b = \eit \wf 1 \dotsb \wf r \cdot 1 \in \binf \cup \ze$.
\end{proof}

\begin{prop}[$\mathbf{O}(r)$]\label{prop:o(r)}
Let $\la \in P^+$ and $\al \in Q_+(r)$.
If $b \in \binf_{-\al}$ satisfies $\ovl \pi_\la (b) \ne 0$, then we have
$\eit \ovl \pi_\la (b) = \ovl \pi_\la (\eit b)$.
\end{prop}
\begin{proof}
Let $P = \wf 1 \dotsb \wf r \cdot 1 \in \linf_{-\al}$,
$b = \wf 1 \dotsb \wf r \cdot 1 + q \linf \in \binf_{-\al}$
and $u = \wf 1 \dotsb \wf r v_\la$.
Then by Proposition \ref{prop:m(r)}, we have
$u \equiv P v_\la \mod q \lla$.

Since $\ovl \pi_\la (b) \ne 0$,
we have $u \notin q \lla$.
Observe that for any $\mu \in P^+$, Lemma \ref{lem:long}(g) gives
\[ \wf 1 \dotsb \wf r (v_\la \ot v_\mu)
\equiv u \ot v_\mu \mod q \lom. \]
Hence by Lemma \ref{lem:eit},
\begin{align}\label{eq:o(r)}
\eit\wf 1 \dotsb \wf r (v_\la \ot v_\mu)
\equiv \eit u \ot v_\mu \mod q \lom.
\end{align}

On the other hand, for any $\mu \gg 0$,
Lemma \ref{lem:e-inf-la} yields
\begin{align*}
\eit (\wf 1 \dotsb \wf r v_{\la+\mu})
\equiv \eit (P v_{\la+\mu})
\equiv (\eit P) v_{\la+\mu} \mod q \llm .
\end{align*}
Applying $\ophi$, Proposition \ref{prop:e(r)} yields
\[ \eit (\wf 1 \dotsb \wf r (v_\la \ot v_\mu))
\equiv (\eit P) (v_\la \ot v_\mu) \mod q \lom . \]
Together with \eqref{eq:o(r)},
we obtain
\[ \eit u\ot v_\mu\equiv (\eit P) (v_\la \ot v_\mu) \mod q \lom. \]
Therefore, by applying $\os$, we conclude
$\eit u \equiv (\eit P) v_\la$, which gives
\[ \eit \ovl \pi_\la (b) = \ovl \pi_\la (\eit b), \]
as claimed.
\end{proof}

\begin{prop}[$\mathbf{B}(r)$]\label{prop:b(r)}
For every $\la \in P^+$ and $\al \in Q_+(r)$, we have
\[ \eit \bla_{\la-\al} \subset \bla \cup \ze. \]
\end{prop}
\begin{proof}
Our assertion is an immediate consequence of
Proposition \ref{prop:o(r)},
Corollary \ref{coro:large-isom}(b) and Proposition \ref{prop:i(r)}.
\end{proof}

To prove the remaining inductive statements,
we use the symmetric form on $\vla$.
Since the situation is special when $a_{ii}=0$,
 we shall introduce
the operator $Q_i$ on $\vla$ by
\begin{align}
Q_i(f_i^{(n)}u)= \begin{cases}
(n+1)f_i^{(n)}u &\mbox{if } a_{ii}=0, \\
f_i^{(n)}u &\mbox{otherwise} \end{cases}
\end{align}
for $u \in \ker e_i$ so that
$Q_i=\id_{\vla}$ when $a_{ii} \ne 0$.

\begin{lem}\label{lem:lla}
Let $\la \in P^+$ be a dominant integral weight
and let $(\ ,\ )$ be the symmetric bilinear form on $\vla$
given by \eqref{eq:bilinear form}. Then we have
\[ (\lla_{\la-\al}, \lla_{\la-\al}) \subset \Az \]
for every $\al \in Q_+(r)$.
Moreover, we have
\begin{align}\label{eq:adjoint}
(\fit u, v) \equiv (u, Q_i\eit v) \mod q \Az
\end{align}
for every $\al \in Q_+(r)$, $u \in \lla_{\la-\al+\al_i}$
and $v \in \lla_{\la-\al}$.
\end{lem}
\begin{proof}
We shall argue by the induction on $|\al|$.
By Corollary \ref{coro:long} (a),
it is enough to show \eqref{eq:adjoint} for
$u = f_i^{(n)} u_0 \neq 0$, $v = f_i^{(m)} v_0 \neq 0$ with
$u_0 \in \ker e_i \cap \lla$ and
$v_0 \in \ker e_i \cap \lmu$.

Set $a=\lan h_i,\wt(v_0)\ran$. If $a=0$ or $m=0$,
then the assertion is trivial. Assume that $a>0$ and $m>0$.
Then we have
\begin{align*}
(\fit u, v)
&= (f_i^{(n+1)} u_0, f_i^{(m)} v_0)
= a (f_i^{(n)} u_0, q_i^{-1}K_i e_i f_i^{(m)} v_0).
\end{align*}
Here, $a=[n+1]_i^{-1}$ or $1$ according as
$a_{ii}=2$ or $a_{ii}\le 0$.
By \eqref{eq:qKe}, we have
\begin{align*}
a q_i^{-1}K_i e_i f_i^{(m)} v_0
\in \begin{cases}
(1+q\Az)f_i^{(m-1)}v_0 &\mbox{if } a_{ii} \ne 0,\\
m(1+q\Az)f_i^{(m-1)}v_0 &\mbox{if } a_{ii}=0.
\end{cases}
\end{align*}
Hence we have
$(\fit u, v)\in(1+q\Az)(u, Q_i\eit v)$.
By the induction hypothesis, we have
$(u, Q_i\eit v)\in\Az$, and we are done.
\end{proof}

\smallskip

\begin{lem}\label{lem:(P,Q)}
Let $\al=\sum_{k=1}^{r}\al_{i_k}\in Q_+$,
$P$, $Q \in \Um_{-\al}$ and $m\in\Z$.
Then for $\la \gg 0$, we have
\[ (P,Q) \equiv
\prod_{k=1}^r (1-q_{i_k}^2)^{-1}(P v_\la, Q v_\la) \mod q^m\Az. \]
\end{lem}
\begin{proof}
If $P = 1$, our assertion is clear.
We shall argue by the induction on the height $r$ of $\al$.
We may assume $P = f_i R$.
Then we have
{\allowdisplaybreaks
\begin{align*}
(f_i R v_\la, Q v_\la)
&= (R v_\la, q_i^{-1} K_i e_i Q v_\la) \\
&= \left( R v_\la, q_i^{-1} K_i \left( Q e_i
+ \frac {K_i^{-1} (e_i' Q)-K_i (e_i'' Q)}{q_i^{-1}-q_i} \right)
v_\la \right) \\
&= \left( R v_\la,
\frac {(e_i' Q) v_\la-K_i^2 (e_i'' Q) v_\la }{1-q_i^2} \right) \\
&= \left( R v_\la,
\frac{(e_i' Q)v_\la - q_i^{2 \lan h_i, \la-\al+\al_i \ran}(e_i'' Q) v_\la}
{1-q_i^2} \right) \\
&=(1-q_i^2)^{-1} (R v_\la, (e_i' Q) v_\la)
-q_i^{2 \lan h_i, \la-\al+\al_i \ran}(1-q_i^2)^{-1}
(R v_\la, (e_i'' Q) v_\la),\\
&\equiv (1-q_i^2)^{-1} (R, e_i' Q)
= (1-q_i^2)^{-1} (f_i R, Q) \mod q^m \Az
\end{align*}
}
by the induction hypothesis.
\end{proof}

For a finitely generated $\Az$-submodule $L$ of
$\vla_{\la-\al}$, we set
\[ L^\vee \seteq \set{u \in \vla_{\la-\al}}{(u, L) \subset \Az}. \]
Then we have $(L^\vee)^\vee=L$.
Similarly, we define $L^\vee$
for a finitely generated $\Az$-submodule $L$ of $\Um_{-\al}$.

\begin{lem}\label{lem:dual}
If $\la \gg 0$, then for any $\al \in Q_+(r)$
we have $\pi_\la (\linf_{-\al}^\vee) = \lla_{\la-\al}^\vee$.
\end{lem}
\begin{proof}
Let $\{ P_k \}$ be an $\Az$-basis of $\linf_{-\al}$,
and $\{ Q_k \}$ the dual basis: $(P_k, Q_j)=\del_{kj}$.
Then $\linf_{-\al}^\vee=\sum_{j} \Az Q_j$.
By Proposition \ref{prop:l(r)}, we have
$\lla_{\la-\al}=\sum_{k} \Az P_k v_\la$.
By Lemma \ref{lem:(P,Q)}, we have
$(P_k v_\la, Q_j v_\la) \equiv \del_{kj} \mod q\Az$ for $\la \gg 0$.
Hence we conclude
\[ \lla_{\la-\al}^\vee=\sum_{j}\Az Q_jv_\la
=\pi_\la(\linf_{-\al}^\vee) \quad\mbox{for }\la\gg 0. \]
\end{proof}

In the following lemma,
we will prove a special case of $\mathbf{F}(r)$.
\begin{lem}\label{lem:special-f(r)}
Let $\mu \gg 0$ and $\al \in Q_+(r)$.
If $\la \in P^+$, then we have
\[ \opsi ((\lom)_{\la+\mu-\al})\subset\llm_{\la+\mu-\al}. \]
\end{lem}
\begin{proof}
Recall that Lemma \ref{lem:phi-1} gives
\[ (\lom)_{\la+\mu-\al}
= \sum_i \fit ((\lom)_{\la+\mu-\al+\al_i})
+ v_\la \ot \lmu_{\mu-\al}. \]
By the induction hypothesis $\mathbf{F}(r-1)$, we have
\begin{align*}
&\opsi ( \sum_i \fit (\lom)_{\la+\mu-\al+\al_i})
= \sum_i \fit \opsi ((\lom)_{\la+\mu-\al+\al_i})\\
&\quad\subset \sum_i \fit \llm_{\la+\mu-\al+\al_i}
= \llm_{\la+\mu-\al}.
\end{align*}

Thus it remains to show
\[ \opsi (v_\la \ot \lmu_{\mu-\al}) \subset \llm_{\la+\mu-\al}. \]
For $u \in \llm_{\la+\mu-\al}^\vee$,
by Lemma \ref{lem:dual}, we can write $u = P v_{\la+\mu}$
for some $P \in \linf_{-\al}^\vee$.
Observe that
\[ \Del(P) = P \ot 1 + (\mbox{intermediate terms}) + K_\al \ot P, \]
which yields
\begin{align*}
(\ophi (u), v_\la \ot \lmu_{\mu-\al})
& = (\Del(P) (v_\la \ot v_\mu), v_\la \ot \lmu_{\mu-\al}) \\
&= (P v_\la \ot v_\mu + \dotsb + K_\al v_\la \ot P v_\mu, v_\la \ot
\lmu_{\mu-\al}) \\
&= q^{(\al| \la)} (P v_\mu, \lmu_{\mu-\al}).
\end{align*}
Since $\mu \gg 0$,
Lemma \ref{lem:dual} implies that $Pv_\mu\in \lmu^\vee$.
Thus we have
\begin{align*}
(u, \opsi (v_\la \ot \lmu_{\mu-\al}))
&= (\ophi (u), v_\la \ot \lmu_{\mu-\al}) \\
&= q^{(\al|\la)} (P v_\mu, \lmu_{\mu-\al}) \subset \Az .
\end{align*}
Hence, we have
\[ \opsi (v_\la \ot \lmu_{\mu-\al})
\subset \llm_{\la+\mu-\al}^{\vee \vee} = \llm_{\la+\mu-\al}, \]
which completes the proof.
\end{proof}

\begin{prop}[$\mathbf{J}(r)$]\label{prop:j(r)}
Suppose $\al \in Q_+(r)$ and $b \in \binf_{-\al}$.
If $\eit b \ne 0$, then we have $b = \fit \eit b$.
\end{prop}

\begin{proof}
Let $b = \wf 1 \dotsb \wf r \cdot 1 \in \binf_{-\al}$
and suppose $\eit b = b' \ne 0$.
If $i_1 = \dotsb = i_r$,
then $\eit b' \ne 0$ implies that
$i_1 = \dotsb = i_r = i$,
and our assertion follows immediately.

If there exists a positive integer $t < r$ such that
$i_t \ne i_{t+1} = \dotsb = i_r$,
take $\mu\gg 0$ and set
$\la_0 = \La_{i_t}$, $\la=\la_0+\mu$.
Then Lemma \ref{lem:key-lemma} yields
\[ \wf 1 \dotsb \wf r (v_{\la_0} \ot v_\mu)
\equiv v \ot v' \mod q L(\la_0) \ot \lmu \]
for some $v \in L(\la_0)_{\la_0-\be}$ and $v' \in \lmu_{\mu-\ga}$
and $\be, \ga \in Q_+(r-1) \setminus \ze$ with $\al = \be + \ga$
such that $v + q L(\la_0) \in B(\la_0) \cup \ze$ and
$v' + q \lmu \in \bmu \cup \ze$.
By Corollary \ref{coro:long} (b), we have
\[ \eit \wf 1 \dotsb \wf r (v_{\la_0} \ot v_\mu)
\equiv \eit (v \ot v') \mod q L(\la_0) \ot \lmu . \]
Applying $\Psi_{\la_0, \mu}$, $\mathbf{G}(r-1)$ yields
\[ \pi_\la(\eit \wf 1 \dotsb \wf r \cdot 1)
= \eit \wf 1 \dotsb \wf r v_{\la_0 + \mu}
\equiv \Psi_{\la_0, \mu}\bigl(\eit(v \ot v')\bigr) \mod q \lla. \]
Thus $\eit (v \ot v') \notin q L(\la_0) \ot \lmu$,
and Lemma \ref{lem:long} (d) gives
\begin{align*}
\wf 1 \dotsb \wf r (v_{\la_0} \ot v_\mu)
&\equiv v \ot v' \equiv \fit \eit (v \ot v') \\
&\equiv \fit \eit \wf 1 \dotsb \wf r (v_{\la_0} \ot v_\mu)
\mod q L(\la_0) \ot \lmu .
\end{align*}
Hence, by applying $\Psi_{\la_0,\mu}$, Lemma \ref{lem:special-f(r)} gives
\[ \wf 1 \dotsb \wf r v_\la
\equiv \fit \eit \wf 1 \dotsb \wf r v_\la \mod q \lla. \]
Therefore by Lemma \ref{lem:infty-la}
and Corollary \ref{coro:large-isom} (c),
we have $b \equiv \fit \eit b \mod q \linf$.
\end{proof}

\begin{prop}[$\mathbf{C}(r)$]\label{prop:c(r)}
Suppose $\al \in Q_+(r)$ and $\la \in P^+$.
Then for every $i \in I$, $b \in \bla_{\la-\al+\al_i}$
and $b' \in \bla_{\la-\al}$,
we have $\fit b = b'$ if and only if $b = \eit b'$.
\end{prop}
\begin{proof}
$(\Longrightarrow)$
Let $\fit b = b' \in \bla_{\la-\al}$.
Since $b \in \bla_{\la-\al+\al_i}$,
by Lemma \ref{lem:la-decomp}(b),
there exists $n \ge 0$ such that
$b \equiv f_i^{(n)} u$ with $e_i u = 0$.
Hence $b' \equiv f_i^{(n+1)} u$ and
$\eit b' \equiv f_i^{(n)} u \equiv b$.

$(\Longleftarrow)$
Let $b' \in \bla_{\la-\al}$ and $b=\eit b'\in\bla$
Then we have $b' = \ovl \pi_\la (b'')$
for some $b'' \in \linf_{-\al}$
by Corollary \ref{coro:large-isom}(b).
By Proposition \ref{prop:o(r)}, we have
\[ \ovl \pi_\la (\eit b'') = \eit \ovl \pi_\la (b'')
= \eit b' \ne 0 . \]
Hence we have $\eit b'' \ne 0$
and Proposition \ref{prop:j(r)} yields
$b'' = \fit \eit b''$.
By applying $\pi_\la$, we obtain
\begin{align*}
\fit \eit b' &= \fit \ovl \pi_\la (\eit b'')
= \ovl \pi_\la (\fit \eit b'') = \ovl \pi_\la (b'') = b',
\end{align*}
as claimed.
\end{proof}

\begin{prop}[$\mathbf{D}(r)$, $\mathbf{K}(r)$]\label{prop:d(r), k(r)}
Let $\al \in Q_+(r)$.
\begin{enumerate}
\item For $\la \in P^+$, $\bla_{\la-\al}$ is a $\Q$-basis
of $\lla_{\la-\al}/q \lla_{\la-\al}$.
\item $\binf_{-\al}$ is a $\Q$-basis of $\linf_{-\al}/q \linf_{-\al}$.
\end{enumerate}
\end{prop}
\begin{proof}
We will prove (a) only.
The proof of (b) is similar.

Suppose that we have a $\Q$-linear dependence relation
\[ \sum_{b \in \bla_{\la-\al}} a_b b = 0 \mbox{ with } a_b \in \Q. \]
Since $\eit \bla_{\la-\al} \subset \bla \cup \ze$ for every $i \in I$, we
get
\[ 0 = \eit \Bigl( \sum_{b \in \bla_{\la-\al}} a_b b \Bigr)
= \sum_{\substack{b \in \bla_{\la-\al},\\ \eit b \ne 0}} a_b (\eit b). \]
By $\mathbf{D}(r-1)$ and $\mathbf{C}(r)$,
the family $\set{\eit b}{b \in \bla_{\la-\al},\ \eit b \ne 0}$ is
linearly independent over $\Q$.
Therefore $a_b = 0$ whenever $\eit b \ne 0$.

But for each $b \in \bla_{\la-\al}$,
there exists $i \in I$ such that $\eit b \ne 0$.
Hence $a_b = 0$ for every $b \in \bla_{\la-\al}$,
which completes the proof.
\end{proof}

\begin{lem}\label{lem:maximal}
Let $\la \in P^+$ and $\al \in Q_+(r)\setminus \ze$.
\begin{enumerate}
\item If $u \in \lla_{\la-\al}/q\lla_{\la-\al}$
satisfies $\eit u = 0$ for every $i \in I$, then $u = 0$.
\item If $u \in \vla_{\la-\al}$ satisfies
$\eit u \in \lla$ for every $i \in I$,
then $u \in \lla$.
\item If $u \in \linf_{-\al}/q \linf_{-\al}$
satisfies $\eit u = 0$ for every $i \in I$, then $u = 0$.
\item If $u \in \Um_{-\al}$ satisfies $\eit u \in \linf$ for every $i \in I$,
then $u \in \linf$.
\end{enumerate}
\end{lem}
\begin{proof}
(a) Write $u = \sum_{b \in \bla_{\la-\al}} a_b b$ with $a_b \in \Q$.
Then for every $i \in I$, we have
\[ \sum_{\eit b \ne 0} a_b (\eit b) = 0. \]
By the same argument as in the proof of Proposition \ref{prop:d(r), k(r)},
all $a_b$ vanish, which implies $u = 0$.

(b) Choose the smallest $N \ge 0$ such that $q^N u \in \lla$.
If $N > 0$, then $\eit (q^N u) = q^N \eit u \in q \lla$ for every $i \in
I$.
Hence, by (a), we would have $q^N u \in q \lla$;
that is, $q^{N-1} u \in \lla$, a contradiction to the minimality of $N$.
Therefore $N = 0$ and $u \in \lla$.

The proofs of (c) and (d) are similar.
\end{proof}

\begin{prop}[$\mathbf{F}(r)$]\label{prop:f(r)}
For every $\la, \mu \in P^+$, we have
\[ \opsi ((\lom)_{\la+\mu-\al}) \subset \llm. \]
\end{prop}
\begin{proof}
We may assume that $\al \ne 0$.
By Corollary \ref{coro:long}(b), we have
\[ \eit ((\lom)_{\la+\mu-\al}) \subset (\lom)_{\la+\mu-\al+\al_i} \]
for every $i \in I$.
Then by applying $\opsi$, the statement $\mathbf{F}(r-1)$ yields
\[ \eit \opsi ((\lom)_{\la+\mu-\al})
\subset \opsi ((\lom)_{\la+\mu-\al+\al_i}) \subset \llm \]
for every $i \in I$.
Hence our assertion follows from Lemma \ref{lem:maximal}(b).
\end{proof}

\begin{prop}[$\mathbf{N}(r)$]\label{prop:n(r)}
For all $\al \in Q_+(r)$ and $\la \in P^+$, we have
\[ \set{b \in \binf_{-\al}}{\ovl \pi_\la (b) \ne 0}
\isoto \bla_{\la-\al}. \]
\end{prop}
\begin{proof}
We may assume $\al \ne 0$.
We already know
$\ovl \pi_\la \binf_{-\al} \setminus \ze = \bla_{\la-\al}$.
Hence it remains to prove that, for $b, b' \in \binf_{-\al}$,
$\ovl \pi_\la (b) = \ovl \pi_\la (b') \ne 0$ implies $b = b'$.
Choose $i \in I$ such that $\eit \ovl \pi_\la b \ne 0$.
By Proposition \ref{prop:o(r)}, we have
\[ \ovl \pi_\la (\eit b) = \ovl \pi_\la (\eit b') \ne 0 . \]
Then by the induction hypothesis $\mathbf{N}(r-1)$ and
Proposition \ref{prop:b(r)},
we have $\eit b = \eit b' \ne 0$.
Hence Proposition \ref{prop:j(r)} yields $b = b'$.
\end{proof}

So far, we have proved all the statements except $\mathbf{G}(r)$.
Using these statements,
we can prove that Lemma \ref{lem:long} holds for all $\al \in Q_+(r)$.
In particular, we have

\begin{lem}\label{lem:revised}\hfill
\begin{enumerate}
\item For all $i \in I$, $\la, \mu \in P^+$ and $\al \in Q_+(r)$, we have
\begin{align*}
\eit (\bla \ot \bmu)_{\mu-\al} \subset (\bla \ot \bmu) \cup \ze.
\end{align*}
\item If $b \ot b' \in (\bla \ot \bmu)_{\mu-\al}$ and $\eit (b \ot b') \ne
0$,
then
\[ b \ot b' \equiv \fit \eit (b \ot b') . \]
\end{enumerate}
\end{lem}

Finally, we complete the grand-loop argument by proving the statement
$\mathbf{G}(r)$.

\begin{prop}[$\mathbf{G}(r)$]\label{prop:g(r)}
For all $i \in I$, $\la, \mu \in P^+$ and $\al \in Q_+(r)$, we have
\[ \opsi ((\bla \ot \bmu)_{\la+\mu-\al}) \subset \blm \cup \ze. \]
\end{prop}
\begin{proof}
We may assume $\al \ne 0$.
Let $b \ot b' \in (\bla \ot \bmu)_{\la+\mu-\al}$.
If $\eit (b \ot b') = 0$ for every $i \in I$,
then $\eit \opsi (b \ot b') = 0$ for every $i \in I$, and
Lemma \ref{lem:maximal} implies $\opsi (b \ot b') = 0$.

If $\eit (b \ot b') \ne 0$ for some $i \in I$,
then by Lemma \ref{lem:revised} and $\mathbf{G}(r-1)$, we have
\begin{align*}
\opsi (b \ot b') &\equiv \opsi (\fit \eit (b \ot b'))
\equiv \fit \opsi (\eit (b \ot b'))\\
&\in \fit (\blm \cup \ze) \subset \blm \cup \ze.
\end{align*}
\end{proof}

Thus the proof of Theorem \ref{thm:existence} is completed.

\bigskip
Let $(\ ,\ )_0$ be the $\Q$-valued symmetric bilinear form
on $\lla / q \lla$ given by
taking the crystal limit of $(\ ,\ )$ on $\lla$. Then we have:
\begin{coro}\label{coro:orthonormal}
Let $\vla$ be the irreducible highest weight $\Uq$-module with highest
weight
$\la \in P^+$
and let $(\lla, \bla)$ be the crystal basis of $\vla$.
\begin{enumerate}
\item The crystal $\bla$ forms an orthogonal basis of $\lla / q \lla$
with respect to $(\ ,\ )_0$,
and $(b,b)_0 \in \Z_{>0}$ for every $b\in\bla$.
In particular, $(\ ,\ )_0$ is positive definite on $\lla / q \lla$.
\item We have
\begin{align*}
\lla &= \set{u \in \vla}{(u, \lla) \subset \Az} \\
&= \set{u \in \vla}{(u, u) \in \Az}.
\end{align*}
\end{enumerate}
\end{coro}
\begin{proof}
(a) We will prove $(b, b')_0\in \del_{b, b'}\Z_{>0}$ for all $b$,
$b' \in \bla_{\la-\al}$,
$\al \in Q_+(r)$ by induction on $|\al|$.
If $|\al|=0$, our assertion is trivial.

If $|\al| > 0$, choose $i \in I$ such that $\eit b \ne 0$.
Then by Theorem \ref{thm:existence} and Proposition \ref{prop:c(r)},
we have
\begin{align*}
(b, b')_0 = (\fit \eit b, b')_0 = (\eit b, Q_i\eit b')_0
\in \Z_{>0}(\eit b,\eit b')_0
\subset\del_{\eit b, \eit b'}\Z_{>0} = \del_{b, b'}\Z_{>0},
\end{align*}
which completes the proof.

(b) Let
\begin{align*}
L_1 &= \set{u \in \vla}{(u, \lla) \subset \Az} , \\
L_2 &= \set{u \in \vla}{(u, u) \in \Az}.
\end{align*}
It is clear that $\lla \subset L_1$ and $\lla \subset L_2$.
Let $u \in \vla$ be such that $(u, \lla) \subset \Az$.
Take the smallest $n \ge 0$ such that $q^n u \in \lla$.
If $n > 0$, then $(q^n u, \lla) \equiv 0 \mod q \Az$.
Since $(\ ,\ )_0$ is non-degenerate, $q^n u \equiv 0 \mod q \lla$;
that is $q^{n-1} u \in \lla$, which contradicts the minimality of $n$.
Hence $n = 0$ and $u \in \lla$, which proves $L_1 \subset \lla$.

Similarly, let $u \in \vla$ be such that $(u, u) \in \Az$.
Take the smallest $n \ge 0$ such that $q^n u \in \lla$.
If $n > 0$, then $(q^n u, q^n u) = q^{2n} (u, u) \in q \Az$.
Since $(\ ,\ )_0$ is positive definite,
we obtain $q^n u \equiv 0 \mod q \lla$,
a contradiction.
Hence $n = 0$ and $u \in \lla$, which proves $L_2 \subset \lla$.
\end{proof}

\begin{rem}
If $a_{ii} \ne 0$ for every $i \in I$, then
$\bla$ is an orthonormal basis of $\lla /q \lla$.
\end{rem}
The following lemma immediately follows from the preceding proposition and
Lemma \ref{lem:(P,Q)}.
\begin{lem}\label{lem:linf}\hfill
\begin{enumerate}
\item $(\linf, \linf) \subset \Az$.
\item Let $(\ ,\ )_0$ be the $\Q$-valued inner product on $\linf/q\linf$
induced by taking the crystal limit of
$(\ , \ )$ on $\linf$.
Then $\binf$ is an orthogonal basis of $\linf/q\linf$.
In particular, $(\ ,\ )_0$ is positive definite on $\linf/q\linf$.
\item We have
\begin{align*}
\linf &= \set{P \in \Um}{(u, \linf) \subset \Az} \\
&= \set{u \in \Um}{(u, u) \in \Az}.
\end{align*}
\end{enumerate}
\end{lem}

\begin{coro}\label{coro:star invariance}
For every $P \in \linf$, we have $P^\star \in \linf$.
\end{coro}
\begin{proof}
Let $P \in \linf$.
Then by Lemma \ref{lem:linf} (a) and \eqref{eq:star},
we have $(P^\star, P^\star) = (P, P) \in \Az$.
Therefore, by Lemma \ref{lem:linf} (c),
we have $P^\star \in \linf$.
\end{proof}

\medskip

\section{{\bf Balanced Triple}}\label{sec:balance}

In the following three sections,
we will {\it globalize} the theory of crystal bases.
Recall that $\Az$ is the subring of $\Q(q)$
consisting of rational functions in $q$
that are regular at $q = 0$.
Similarly, let $\Ai$ be the subring of $\Q(q)$
consisting of rational functions in $q$
that are regular at $q = \infty$.
Thus we have
\[ \Az / q \Az \isoto \Q \quad \mbox{and}
\quad \Ai / q^{-1} \Ai \isoto \Q \]
by evaluation at $q = 0$ and $q = \infty$, respectively.
Finally, we denote by $\A = \Q[q, q^{-1}]$,
the ring of Laurent polynomials in $q$.

For a subring $A$ of $\Q(q)$ and a vector space $V$ over $\Q(q)$,
an $A$-{\em lattice} is, by definition,
a free $A$-module of $V$ generating $V$ as a vector space over $\Q(q)$.
Assume that
\begin{align}\label{eq:VA}
\parbox{350pt}
{$V$ is a vector space over $\Q(q)$, and $V^\A$, $L_0$ and $\Li$ are
an $\A$-lattice, an $\Az$-lattice and an $\Ai$-lattice of $V$,
respectively.}
\end{align}
Note that $(V^\A\cap L_0)/(V^\A\cap qL_0)\to L_0/qL_0$ is an isomorphism.

\begin{defi}
Set $E = V^\A \cap L_0 \cap \Li$.
Then a triple $(V^\A, L_0, \Li)$ is called a {\it balanced triple} for $V$
if the canonical homomorphisms
$\A\ot E\to V^\A$,
$\Az\ot E\to L_0$ and
$\Ai\ot E\to \Li$ are isomorphisms.
\end{defi}

\begin{thm}\label{thm:equivalence VT}
Let $V$, $V^\A$, $L_0$ and $\Li$ be as in \eqref{eq:VA}.
Set $E = V^\A \cap L_0 \cap \Li$.
Then the following statements are equivalent.
\begin{enumerate}
\item $(V^\A, L_0, \Li)$ is a balanced triple.
\item The canonical map $E \to L_0 / q L_0$ is an isomorphism.
\item The canonical map $E \to \Li / q^{-1} \Li$ is an isomorphism.
\end{enumerate}
\end{thm}

Let $(V^\A, L_0, \Li)$ be a balanced triple for a $\Q(q)$-vector space $V$,
and let
\[ G \cl L_0 / q L_0 \longrightarrow E \seteq V^\A \cap L_0 \cap \Li \]
be the inverse of the canonical isomorphism $E \isoto L_0 / q L_0$.
If $B$ is a $\Q$-basis of $L_0 / q L_0$,
then $G(B) \seteq \set{G(b)}{b \in B }$ is a $\Q$-basis of $E$.
Therefore $G(B)$ is a $\Q(q)$-basis of $V$, an $\A$-basis of $V^\A$,
an $\Az$-basis of $L_0$ and an $\Ai$-basis of $\Li$ as well.

\begin{defi}
We call $G(B)$ the {\it global basis} of $V$
corresponding to the {\it local basis} $B$ of $V$ at $q = 0$.
\end{defi}

Our goal is to construct the global basis $\gla$ of $\vla$
(resp. $\ginf$ of $\Um$)
corresponding to the local basis $\bla$ of $\vla$
(resp. $\binf$ of $\Um$) at $q = 0$.
For this purpose, we need the following lemma which is proved in
\cite{Kas91}.

\begin{lem}[\cite{Kas91}]\label{lem:balanced}
Let $V$, $V^\A$, $L_0$ and $\Li$ be as in \eqref{eq:VA}.
Let $F$ be a $\Q$-vector subspace of $V^\A\cap L_0 \cap \Li$
satisfying the following conditions{\rm :}
\begin{enumerate}
\item the canonical maps
$F \longrightarrow L_0 / q L_0$ and
$F \longrightarrow \Li / q^{-1} \Li$
are injective.
\item $V^\A = \A F$.
\end{enumerate}
Then $(V^\A,L_0,\Li)$ is a balanced triple and $F=V^\A\cap L_0 \cap \Li$.
\end{lem}

\medskip

\section{{\bf Global Bases}}

Let $U_\A^0 (\g)$ be the $\A$-subalgebra of $\Uq$ generated by
$q^h$, $\prod_{k=1}^m \dfrac {1-q^k\cdot q^h}{1-q^k}$
($h \in \pv$, $m\in\Z_{>0}$).
We denote by $U_\A^+ (\g)$ (resp.\ $U_\A^- (\g)$)
the $\A$-subalgebra of $\Uq$
generated by $e_i^{(n)}$ (resp.\ $f_i^{(n)}$)
for $i \in I$, $n \in \Z_{\ge 0}$.
Let $U_\A (\g)$ be the $\A$-subalgebra of $\Uq$ generated by
$U_\A^0 (\g)$, $U_\A^+ (\g)$ and $U_\A^- (\g)$.

Then we can show that
\begin{align}\label{eq:Um decomp}
U_\A (\g) \simeq U_\A^- (\g) \ot U_\A^0 (\g) \ot U_\A^+ (\g) .
\end{align}

Let $\vla$ be the irreducible highest weight
$\Uq$-module with highest weight $\la \in P^+$
and highest weight vector $v_\la$, and let
\[ \vla^\A = U_\A (\g) v_\la = U_\A^- (\g) v_\la . \]

Consider the $\Q$-algebra automorphism ${ }^- \cl \Uq \to \Uq$ defined by
\begin{align}\label{eq:auto}
q \mt q^{-1}, \quad q^h \mt q^{-h}, \quad e_i \mt e_i, \quad f_i \mt f_i
\end{align}
for $h \in \pv$ and $i \in I$.
Then we get a $\Q$-linear automorphism ${ }^-$ on $\vla$ defined by
\begin{align}\label{eq:auto-1}
P v_\la \mt \ovl P v_\la \quad \mbox{for} \quad P \in \Uq.
\end{align}
Our goal is to prove the following theorem.

\begin{thm}\label{thm:balanced triple}\hfill
\begin{enumerate}
\item $(U_\A^- (\g), \linf, \linf^-)$ is a balanced triple for $\Um$.
\item $(\vla^\A, \lla, \lla^-)$ is a balanced triple for $\vla$.
\end{enumerate}
\end{thm}
Once we have Theorem \ref{thm:balanced triple}, the canonical maps
\begin{align*}
&U_\A^-(\g) \cap \linf \cap \linf^- \longrightarrow \linf / q \linf , \\
&\vla^\A \cap \lla \cap \lla^- \longrightarrow \lla / q \lla
\end{align*}
are isomorphisms.
As in Section \ref{sec:balance}, let
\begin{align*}
G_\infty&\cl \linf / q \linf \longrightarrow
U_\A^-(\g) \cap \linf \cap \linf^- ,\\
G_\la&\cl \lla / q \lla \longrightarrow
\vla^\A \cap \lla \cap \lla^-
\end{align*}
be the inverses of the isomorphisms above, and set
\begin{align} \label{eq:global}
\ginf = \set{G_\infty (b)}{b \in \binf}, \quad
\gla = \set{G_\la(b)} {b \in \bla}.
\end{align}
Then we have:
\begin{thm}\hfill
\begin{enumerate}
\item $\ginf$ $(\mbox{resp.\ } \gla)$ is
an $\A$-basis of $U_\A (\g)$ $(\mbox{resp.\ } \vla^\A)$.
\item $\ginf$ $(\mbox{resp.\ } \gla)$ is
an $\Az$-basis of $\linf$ $(\mbox{resp.\ } \lla)$.
\item $\ginf$ $(\mbox{resp.\ } \gla)$ is
a $\Q(q)$-basis of $\Um$ $(\mbox{resp.\ } \vla)$.
\item $\ginf$ $(\mbox{resp.\ } \gla)$ is
a $\Q[q]$-basis of $U_\A^- (\g) \cap \linf$
$(\mbox{resp.\ } \vla^\A \cap \lla)$.
\end{enumerate}
\end{thm}
\begin{proof}
These are immediate consequences of the definition of a balanced triple.
\end{proof}
Moreover, we have
\begin{thm}\label{thm:global auto}\hfill
\begin{enumerate}
\item For $b\in\binf$,
$G_\infty(b)$ is a unique element of $U_\A^- (\g)\cap\linf$ such that
\[ G_\infty (b) \equiv b \mod q \linf,
\quad \ovl {G_\infty (b)} = G_\infty (b). \]
\item For $b\in\bla$, $G_\la(b)$ is a unique element of $\vla^\A\cap \lla$
such that
\[ G_\la (b) \equiv b \mod q \lla, \quad \ovl {G_\la (b)} = G_\la (b). \]
\end{enumerate}
\end{thm}
\begin{proof}
The proof of (a) being similar, we will prove only (b).
By the definition of $G_\la$, it is clear that
$G_\la (b) \equiv b\mod q \lla$.
Set
\[ v = \frac{G_\la (b) - \ovl{G_\la (b)}}{q - q^{-1}} . \]
Clearly, $G_\la (b)$ and $\ovl {G_\la (b)}$ are contained in
$E \seteq \vla^\A \cap \lla \cap
\lla^-$.
Write $G_\la (b) = \sum_j c_j (q) P_j v_\la$, where $c_j (q) \in \A$ and
$P_j$ is a monomial in $f_i^{(n)}$'s ($i \in I$, $n \in \Z_{>0}$).
Then we have $\ovl {G_\la (b)} = \sum_j c_j (q^{-1}) P_j v_\la$,
from which it follows that
\[ v \seteq \frac{G_\la (b) - \ovl{G_\la (b)}}{q - q^{-1}}
= \sum_j \frac{c_j(q) - c_j(q^{-1})}{q - q^{-1}} P_j v_\la \in \vla^\A. \]
Since $\dfrac 1 {q - q^{-1}} = \dfrac q {q^2 - 1} \in q \Az \cap \Ai$, we
have
$v\in q \lla \cap \lla^-$.
Thus $v$ is mapped to $0$ under the canonical isomorphism
$ E\to \lla / q \lla$,
which implies $v= 0$. Therefore we get
$\ovl {G_\la (b)} = G_\la (b)$.

Finally, suppose that $u\in\vla^\A\cap\lla$ satisfies
\[ u \equiv b \mod q \lla \quad \mbox{and}\quad \ovl u = u . \]
Then $u\in E \seteq \vla^\A \cap \lla \cap \lla^-$.
Since the canonical map $E \to \lla / q \lla$ is
an isomorphism,
we must have $u = G_\la (b)$.
\end{proof}

\begin{defi}
We call $\ginf$ $(\mbox{resp.\ } \gla)$
the {\it global basis} of $\Um$ $(\mbox{resp.\ } \vla)$.
\end{defi}

\medskip

\section{{\bf Existence of Global Bases}}

In this section, we will prove
Theorem \ref{thm:balanced triple}.
For $i \in I$ and $n \ge 1$, define
\[ (f_i^n \Um)^\A = f_i^n \Um \cap U_\A^- (\g), \quad
(f_i^n \vla)^\A = (f_i^n \Um)^\A\vla. \]
Then $(f_i^n \vla)^\A$ is an $\A$-lattice of $f_i^n \vla$.

\begin{lem}\label{lem:g decomp}
Fix $i \in I$.
Let $u \in \Um$ and consider the \str\
$u = \sum_{l \ge 0} f_i^{(l)} u_l$ with $e_i' u_l = 0$ for every $l \ge 0$.
If $u \in U_\A^- (\g)$, then $u_l \in U_\A^- (\g)$ for every $l \ge 0$.
\end{lem}
\begin{proof}
We first claim that
\begin{align}
e_i'{}^{(n)} U_\A^- (\g)\subset U_\A^- (\g).
\end{align}

\noindent
Indeed, the commutation relation \eqref{eq:come'f}
implies that
\[ \set{Q \in \Um}{e_i'{}^{(n)} Q \in \Um^\A
\mbox{ for every } n \ge 0} \]
is stable under the left multiplication by
any $f_j^{(m)}$ $(j \in I, m \ge 0)$.

Hence the projection operator $P_i$ introduced in \eqref{eq:proj}
satisfies $P_iU_\A^- (\g)\subset U_\A^- (\g)$.
Then Proposition \ref{prop:decomp}
implies that all $u_l$ belong to $U_\A^- (\g)$.
\end{proof}

Then as an immediate consequence of Lemma \ref{lem:g decomp}, we have:
\begin{lem}\hfill
\begin{enumerate}
\item $(f_i^n \Um)^\A \seteq f_i^n \Um \cap U_\A^- (\g)
= \sum_{l \ge n} f_i^{(l)} U_\A^- (\g)$.
\item $(f_i^n \vla)^\A \seteq (f_i^n \Um)^\A\vla
= \sum_{l \ge n} f_i^{(l)} \vla^\A$.
\end{enumerate}
\end{lem}
In order to prove Theorem \ref{thm:balanced triple},
we shall show the following statements by the induction on $r$.

\begin{enumerate}\item[]
\begin{enumerate}
\item[$\mathbf{X}(r):$] For every $\al \in Q_+(r)$,
the triple $(U_\A^- (\g)_{-\al},\linf_{-\al},\linf^-_{-\al})$
is balanced.
\item[$\mathbf{Y}(r):$]
For all $\la \in P^+$ and $\al \in Q_+(r)$,
the triple $(\vla_{\la-\al},\lla_{\la-\al},\lla^-_{\la-\al})$
is balanced.
\item[$\mathbf{Z}(r):$] For every $\al \in Q_+(r)$ and $\la\in P^+$, let
\[ G_\la \cl \lla _{\la-\al}/q \lla_{\la-\al}
\longrightarrow \vla^\A_{\la-\al}\cap\lla\cap\lla^- \]
be the inverse of the canonical isomorphism
\[ \vla^\A_{\la-\al}\cap\lla\cap\lla^- \isoto
\lla _{\la-\al}/q \lla_{\la-\al} \]
guaranteed by $\mathbf{Y}(r)$.
Then for any $b \in \bla_{\la-\al} \cap \fit^n \bla$ with $n \ge 0$,
we have $G_\la(b) \in (f_i^n \vla)^\A$.
\end{enumerate}
\end{enumerate}

If $r=0$, these assertions are obvious.
Suppose $r > 0$ and assume that the statements
$\mathbf{X}(r-1)$, $\mathbf{Y}(r-1)$ and $\mathbf{Z}(r-1)$
are true.
Then we can define
\begin{align*}
G_\infty &\cl \linf_{-\be}/q\linf_{-\be}
\isoto U_\A^-(\g)_{-\be}\cap \linf\cap\linf^-,\\
G_\la &\cl \lla_{\la-\be}/q\lla_{\la-\be}
\isoto \vla^\A_{\la-\be}\cap \lla\cap\lla^-
\end{align*}
for $\be\in Q_+(r-1)$.
Thus we have:
\begin{lem}\label{lem:g basic}
For all $\la \in P^+$ and $\be \in Q_+(r-1)$, we have
\begin{enumerate}
\item $U_\A^- (\g)_{-\be} = \bop_{b \in \binf_{-\be}} \A G_\infty(b)$,\\
$\vla_{\la-\be}^\A = \bop_{b \in \bla_{\la-\be}} \A G_\la(b)$,
\item $\ovl {G_\infty(b)} = G_\infty(b)$
for every $b \in \binf_{-\be}$,\\
$\ovl {G_\la(b)} = G_\la(b)$ for every $b \in \bla_{\la-\be}$,
\item $G_\infty (b) v_\la = G_\la (\ovl \pi_\la (b))$
for every $b \in \binf_{-\be}$.
\end{enumerate}
\end{lem}
The following proposition will play a crucial role in proving the statement
$\mathbf{X}(r)$ and $\mathbf{Y}(r)$.

\begin{prop}\label{prop:g main}
Let $i \in I$, $\al \in Q_+(r)$. Then for any $n \ge 1$,
the triple
\[ \Bigl((f_i^n \vla)_{\la-\al}^\A ,f_i^n \vla\cap\lla_{\la-\al},
f_i^n \vla\cap\lla^-_{\la-\al}\Bigr) \]
is balanced.
\end{prop}
\begin{rem}
By Proposition \ref{prop:fn}, we have
\[ \bigl(f_i^n \vla\cap\lla_{\la-\al}\bigr)
/\bigl(f_i^n \vla\cap q\lla_{\la-\al}\bigr)
\simeq \bop_{b\in \bla_{\la-\al}\cap\fit^n\bla}\Q b. \]
\end{rem}
\begin{proof}
If $a_{ii} = 2$, our assertion can be proved as in \cite{Kas91}.
Hence we assume that $a_{ii} \le 0$.
In this case, we have
\begin{align*}
(f_i^n \vla)_{\la-\al}^\A
= \sum_{l \ge n} f_i^l \bigl(\vla_{\la-\al+l\al_i}^\A\bigr)
= f_i^n\bigl( \vla_{\la-\al+n\al_i}^\A\bigr).
\end{align*}
We will use the descending induction on $n$.
Note that if $n > r$, then $(f_i^n \vla)_{\la-\al}= 0$,
and our assertion is trivial.
Assume that our assertion is true for $n+1$.
Set
\[ V_n^\A = (f_i^{n} \vla)_{\la-\al}^\A,\
L_n = f_i^n \vla \cap \lla_{\la-\al} \mbox{ and }
E_n = V_n^\A \cap L_n\cap L_n^-. \]

If $\lan h_i,\la-\al+n\al_i \ran = 0$, then
$f_i^n \bigl(\vla_{\la-\al+n\al_i}\bigr)=0$,
and our assertion is trivial.
Hence we may assume that $\lan h_i,\la-\al+n\al_i\ran>0$,
in which case, $\fit^n\cl \bla_{\la-\al+n\al_i}\to
\bla_{\la-\al}$ is injective.

Since $n \ge 1$, Lemma \ref{lem:g basic}(a) gives
$\vla_{\la-\al+n\al_i}^\A = \bop_{b \in \bla_{\la-\al+n\al_i}} \A G_\la(b)$
and hence we have
\[ V_n=\sum_{b \in \bla_{\la-\al+n\al_i}} \A f_i^n G_\la(b). \]
If $b \in \bla_{\la-\al+n\al_i}$ satisfies $\eit b \ne 0$, then
$\mathbf{Z}(r-1)$ implies that
$G_\la(b) \in (f_i \vla)^\A = f_i \vla^\A$, which yields
\[ f_i^n G_\la(b) \in V_{n+1}^\A \seteq (f_i^{n+1} \vla)_{\la-\al}^\A . \]

\smallskip

Set $B_0 \seteq \set{b \in \bla_{\la-\al+n\al_i}}{\eit b = 0 }$ and
$F_0 \seteq \sum_{b \in B_0} \Q f_i^n G_\la(b) $.
Then we get
\begin{align}\label{eq:V_n^A}
V_n^\A= \A F_0+ V_{n+1}^\A.
\end{align}
Note that
\begin{align*}
f_i^n G_\la(b)\in E_n, \quad
f_i^n G_\la(b)\equiv \fit^nb\mod q\lla
\quad \mbox{for every } b \in B_0.
\end{align*}
Set $F=F_0+E_{n+1}$.
Then $F\subset E_n$, and $V_{n+1}^\A=\A E_{n+1}$ by the induction
hypothesis.
Hence \eqref{eq:V_n^A} implies
$V_n^\A=\A F$.

Consider the morphism
\[ F \to L_n/qL_n
\simeq \Bigl(\bop_{b\in B_0}\Q\fit^nb\Bigr)\op (L_{n+1}/qL_{n+1}). \]
Since $F_0\to \op_{b\in B_0}\Q\fit^nb$ and
$E_{n+1}\to L_{n+1}/qL_{n+1}$ are isomorphisms,
$F\to L_n/qL_n$ is an isomorphism.
Similarly, $F\to (L_n)^-/q^{-1}(L_n)^-$ is an isomorphism.
Hence Lemma \ref{lem:balanced} implies the desired result.
\end{proof}

\begin{coro}\label{coro:g coro}
For every $i \in I$, $\al \in Q_+(r)$ and $n \ge 1$,
we have the following canonical isomorphism\,{\rm :}
\begin{align*}
(f_i^n \Um)_{-\al}^\A \cap \linf \cap \linf^-
&\isoto \frac {(f_i^n \Um)_{-\al}^\A \cap \linf}
{(f_i^n \Um)_{-\al}^\A \cap q \linf} \\
&\simeq \bop_{b \in \binf_{-\al} \cap \fit^n \binf} \Q b.
\end{align*}
\end{coro}
\begin{proof}
It suffices to observe that, for $\la \gg 0$, we have
\begin{align*}
&\Um_{-\al} \isoto \vla_{\la-\al},
\quad (f_i^n \Um)_{-\al}^\A \isoto (f_i^n \vla)^\A_{\la-\al}, \\
&\linf_{-\al} \isoto \lla_{\la-\al},
\quad \linf_{-\al}^- \isoto \lla_{\la-\al}^-,\\
&\bop_{b \in \binf_{-\al} \cap \fit^n \binf} \Q b
\isoto \bop_{b \in \bla_{\la-\al} \cap \fit^n \bla} \Q b .
\end{align*}
\end{proof}
For $\al \in Q_+(r)$, let
\[ G_i \cl \bop_{b \in \binf_{-\al} \cap \fit \binf} \Q b
\longrightarrow (f_i \Um)_{-\al}^\A \cap \linf \cap \linf^- \]
be the inverse of the canonical isomorphism
given in Corollary \ref{coro:g coro} (with $n = 1$).
Then we have
\[ (f_i \Um)_{-\al}^\A
= \bop_{b \in \binf_{-\al} \cap \fit \binf} \A G_i (b). \]
Now we will prove:
\begin{lem}\label{lem:well global}
For $i$, $j\in I$, $\al \in Q_+(r)$ and
$b \in \binf_{-\al} \cap \fit \binf \cap \tilde f_j \binf$,
we have
\[ G_i (b) = G_j (b) . \]
\end{lem}
\begin{proof}
Let $b = \fit \dotsb \tilde f_k \cdot 1
\in \binf_{-\al} \cap \fit \binf \cap \tilde f_j \binf$.
Take $\la \in P^+$ such that
$\la (h_k) = 0$ and $\la(h_p) \gg 0$ for $p\in I\setminus\{k\}$.
Since $f_kv_\la=0$,
we get $\ovl \pi_\la (b) = 0$,
which implies
\[ G_i (b) v_\la \equiv 0 \mod q \lla. \]
That is, $G_i (b) v_\la \in (f_i \vla)_{\la-\al}^\A \cap q \lla \cap \lla^-$.

By Proposition \ref{prop:g main},
we must have $G_i (b) v_\la = 0$ in $\vla$.
Note that
\begin{align*}
\vla_{\la-\al} &\simeq \Um_{-\al} /\bigl(\Um_{-\al+\al_k} f_k
+\sum_{p \in I^\re \setminus \{ k \}} \Um_{-\al+(\la(h_p)+1) \al_i}
f_p^{(\la(h_p)+1)}\bigl)\\
&\simeq \Um_{-\al} /\Um_{-\al+\al_k} f_k
\end{align*}
as a $\Um$-module.
It follows that
$G_i (b) \in \Um f_k \cap U_\A^- (\g)$.
By applying the anti-automorphism $\star$ on $G_i(b)$, we have
\[ G_i(b)^\star \in f_k \Um\cap U_\A^- (\g)
= (f_k \Uq)^\A . \]
Moreover, by Corollary \ref{coro:star invariance},
we have $G_i(b)^\star \in (f_k \Um)_{-\al}^\A \cap \linf \cap \linf^-$.

Similarly,
$G_j(b)$ satisfies $G_j(b)^\star \in (f_k \Um)^\A \cap \linf \cap \linf^-$.
Since $G_i(b) \equiv G_j(b) \equiv b \mod q \linf$, we have
$G_i(b)^\star\equiv G_j(b)^\star\mod q \linf$.
Hence, Corollary \ref{coro:g coro} implies
$G_i(b)^\star = G_j(b)^\star$, which yields $G_i(b) = G_j(b)$ as desired.
\end{proof}

\smallskip

Let $b \in \binf_{-\al}$ with $\al \in Q_+(r)$.
Since $r \ge 1$, we have $\eit b \ne 0$ for some $i \in I$;
that is, $b \in \binf_{-\al} \cap \fit \binf$ for some $i \in I$.
By Lemma \ref{lem:well global}, we may define a map
\[ G_\infty \cl \linf_{-\al} / q \linf_{-\al}
\longrightarrow \Um_{-\al}^\A \cap \linf \cap \linf^- \]
by $G_\infty (b) = G_i (b)$ for any $i$ with $\eit b \ne 0$.
By the definition, we have
\begin{align}
\begin{aligned}
&G_\infty (b) \equiv b \mod q \linf,\\
&(f_i^n \Um)_{-\al}^\A
= \bop_{b \in \binf_{-\al} \cap \fit^n \binf} \A G_\infty (b)
\quad\mbox{for } n \ge 1 .
\end{aligned}
\end{align}
Set $F = \sum_{b \in \binf_{-\al}} \Q G_\infty(b)
\subset \Um^\A_{-\al}\cap \linf\cap\linf^-$.
Since $U_\A^- (\g)_{-\al} = \sum_i (f_i \Um)_{-\al}^\A$,
we obtain
\[ U_\A^- (\g)_{-\al} = \A F. \]
Moreover, the canonical maps
\[ F \to \linf_{-\al}/q\linf_{-\al}
\mbox{ and } F \to \linf_{-\al}^-/q^{-1}\linf_{-\al}^- \]
are bijective.
Therefore, by Lemma \ref{lem:balanced},
we conclude $\mathbf{X}(r)$.

Let us prove $\mathbf{Y}(r)$.
Note that for all $\al \in Q_+(r)$, $b \in \binf_{-\al}$ and $\la \in P^+$,
if $\ovl \pi_\la (b) = 0$, then $G_\infty (b) v_\la = 0$.
Indeed, let $i \in I$ be such that $\eit b \ne 0$.
Then we have
\[ G_\infty (b) v_\la \in
(f_i \vla)_{\la-\al}^\A \cap q \lla \cap \lla^- , \]
and by Proposition \ref{prop:g main},
we obtain $G_\infty (b) v_\la = 0$.

Set $F \seteq \sum
\limits_{\substack{b \in \binf_{-\al}\\ \ovl \pi_\la(b) \ne 0}}
\Q G_\infty (b) v_\la\subset \vla^\A_{\la-\al}\cap\lla\cap\lla^-$.
Then we have
$\vla_{\la-\al}^\A = \A F$.
By Proposition \ref{prop:n(r)}, the homomorphisms
$F\to \lla_{\la-\al}/q\lla_{\la-\al}$
and
$F\to \lla_{\la-\al}^-/q^{-1}\lla_{\la-\al}^-$
are isomorphisms.
Therefore, by Lemma \ref{lem:balanced},
we obtain the statement $\mathbf{Y}(r)$.

Finally, the statement $\mathbf{Z}(r)$ follows from
Proposition \ref{prop:g main},
which completes the proof.

\begin{rem}\hfill
\begin{tenumerate}
\item We have $G_\infty(b) v_\la=G_\la(\ovl \pi_\la (b))$
for any $\la\in P^+$ and $b\in\binf$.
\item
If $i \in I^\im$, we have
\[ G_\la(\fit^n b)=f_i^nG_\la(b) \mbox{ for every } b \in \bla
\mbox{ and } n \in \Z_{\ge 0}. \]
\item We have
\[ (f_i^n\vla)^\A=f_i^n\vla\cap\vla^\A. \]
Indeed, $\{G_\la(b)\}_{b\in\fit^n \bla}$ is an $\A$-basis of
$(f_i^n\vla)^\A$, and hence it is also a $\Q(q)$-basis of $f_i^n \vla$.
\end{tenumerate}
\end{rem}

\end{document}